\documentclass[11pt]{article}
\usepackage{latexsym,amssymb,amsmath,ulem}
\usepackage{color}
\usepackage{graphicx}
\usepackage{amsfonts}

\def\tto{\;{\lower 1pt \hbox{$\rightarrow$}}\kern -10pt
\hbox{\raise 2pt \hbox{$\rightarrow$}}\;}

\newtheorem{theorem}{Theorem}[section]
\newtheorem{proposition}{Proposition}[section]
\newtheorem{corollary}{Corollary}[section]
\newtheorem{lemma}{Lemma}[section]
\newtheorem{remark}{Remark}[section]
\newtheorem{example}{Example}[section]
\newtheorem{definition}{Definition}[section]

\numberwithin{equation}{section}

\renewcommand{\theequation}{\thesection.\arabic{equation}}
\normalsize

\makeatletter

\ifx\pdfoutput\relax\let\pdfoutput=\undefined\fi
\newcount\msipdfoutput
\ifx\pdfoutput\undefined
\else
 \ifcase\pdfoutput
 \else
    \ifx\paperwidth\undefined
    \else
      \ifdim\paperheight=0pt\relax
      \else
        \pdfpageheight\paperheight
      \fi
      \ifdim\paperwidth=0pt\relax
      \else
        \pdfpagewidth\paperwidth
      \fi
    \fi
  \fi
\fi

\newcount\@hour\newcount\@minute\chardef\@x10\chardef\@xv60
\def\tcitime{
\def\@time{%
  \@minute\time\@hour\@minute\divide\@hour\@xv
  \ifnum\@hour<\@x 0\fi\the\@hour:%
  \multiply\@hour\@xv\advance\@minute-\@hour
  \ifnum\@minute<\@x 0\fi\the\@minute
  }}%

\def\x@hyperref#1#2#3{%
   \catcode`\~ = 12
   \catcode`\$ = 12
   \catcode`\_ = 12
   \catcode`\# = 12
   \catcode`\& = 12
   \catcode`\% = 12
   \y@hyperref{#1}{#2}{#3}%
}

\def\y@hyperref#1#2#3#4{%
   #2\ref{#4}#3
   \catcode`\~ = 13
   \catcode`\$ = 3
   \catcode`\_ = 8
   \catcode`\# = 6
   \catcode`\& = 4
   \catcode`\% = 14
}

\@ifundefined{hyperref}{\let\hyperref\x@hyperref}{}
\@ifundefined{msihyperref}{\let\msihyperref\x@hyperref}{}

\@ifundefined{qExtProgCall}{\def\qExtProgCall#1#2#3#4#5#6{\relax}}{}
\def\QCTOpt[#1]#2{%
  \def\QCTOptB{#1}
  \def\QCTOptA{#2}
}
\def\QCTNOpt#1{%
  \def\QCTOptA{#1}
  \let\QCTOptB\empty
}
\def\Qct{%
  \@ifnextchar[{%
    \QCTOpt}{\QCTNOpt}
}
\def\QCBOpt[#1]#2{%
  \def\QCBOptB{#1}%
  \def\QCBOptA{#2}%
}
\def\QCBNOpt#1{%
  \def\QCBOptA{#1}%
  \let\QCBOptB\empty
}
\def\Qcb{%
  \@ifnextchar[{%
    \QCBOpt}{\QCBNOpt}%
}
\def\PrepCapArgs{%
  \ifx\QCBOptA\empty
    \ifx\QCTOptA\empty
      {}%
    \else
      \ifx\QCTOptB\empty
        {\QCTOptA}%
      \else
        [\QCTOptB]{\QCTOptA}%
      \fi
    \fi
  \else
    \ifx\QCBOptA\empty
      {}%
    \else
      \ifx\QCBOptB\empty
        {\QCBOptA}%
      \else
        [\QCBOptB]{\QCBOptA}%
      \fi
    \fi
  \fi
}
\newcount\GRAPHICSTYPE
\GRAPHICSTYPE=\z@
\def\GRAPHICSPS#1{%
 \ifcase\GRAPHICSTYPE
   \special{ps: #1}%
 \or
   \special{language "PS", include "#1"}%
 \fi
}%
\def\graffile#1#2#3#4{%
    \bgroup
	   \@inlabelfalse
       \leavevmode
       \@ifundefined{bbl@deactivate}{\def~{\string~}}{\activesoff}%
        \raise -#4 \BOXTHEFRAME{%
           \hbox to #2{\raise #3\hbox to #2{\null #1\hfil}}}%
    \egroup
}%
\def\draftbox#1#2#3#4{%
 \leavevmode\raise -#4 \hbox{%
  \frame{\rlap{\protect\tiny #1}\hbox to #2%
   {\vrule height#3 width\z@ depth\z@\hfil}%
  }%
 }%
}%
\newcount\@msidraft
\@msidraft=\z@
\let\nographics=\@msidraft
\newif\ifwasdraft
\wasdraftfalse

\def\GRAPHIC#1#2#3#4#5{%
   \ifnum\@msidraft=\@ne\draftbox{#2}{#3}{#4}{#5}%
   \else\graffile{#1}{#3}{#4}{#5}%
   \fi
}
\def\addtoLaTeXparams#1{%
    \edef\LaTeXparams{\LaTeXparams #1}}%
\newif\ifBoxFrame \BoxFramefalse
\newif\ifOverFrame \OverFramefalse
\newif\ifUnderFrame \UnderFramefalse

\def\BOXTHEFRAME#1{%
   \hbox{%
      \ifBoxFrame
         \frame{#1}%
      \else
         {#1}%
      \fi
   }%
}

\def\doFRAMEparams#1{\BoxFramefalse\OverFramefalse\UnderFramefalse\readFRAMEparams#1\end}%
\def\readFRAMEparams#1{%
 \ifx#1\end%
  \let\next=\relax
  \else
  \ifx#1i\dispkind=\z@\fi
  \ifx#1d\dispkind=\@ne\fi
  \ifx#1f\dispkind=\tw@\fi
  \ifx#1t\addtoLaTeXparams{t}\fi
  \ifx#1b\addtoLaTeXparams{b}\fi
  \ifx#1p\addtoLaTeXparams{p}\fi
  \ifx#1h\addtoLaTeXparams{h}\fi
  \ifx#1X\BoxFrametrue\fi
  \ifx#1O\OverFrametrue\fi
  \ifx#1U\UnderFrametrue\fi
  \ifx#1w
    \ifnum\@msidraft=1\wasdrafttrue\else\wasdraftfalse\fi
    \@msidraft=\@ne
  \fi
  \let\next=\readFRAMEparams
  \fi
 \next
 }%
\def\IFRAME#1#2#3#4#5#6{%
      \bgroup
      \let\QCTOptA\empty
      \let\QCTOptB\empty
      \let\QCBOptA\empty
      \let\QCBOptB\empty
      #6%
      \parindent=0pt
      \leftskip=0pt
      \rightskip=0pt
      \setbox0=\hbox{\QCBOptA}%
      \@tempdima=#1\relax
      \ifOverFrame
          \typeout{This is not implemented yet}%
          \show\HELP
      \else
         \ifdim\wd0>\@tempdima
            \advance\@tempdima by \@tempdima
            \ifdim\wd0 >\@tempdima
               \setbox1 =\vbox{%
                  \unskip\hbox to \@tempdima{\hfill\GRAPHIC{#5}{#4}{#1}{#2}{#3}\hfill}%
                  \unskip\hbox to \@tempdima{\parbox[b]{\@tempdima}{\QCBOptA}}%
               }%
               \wd1=\@tempdima
            \else
               \textwidth=\wd0
               \setbox1 =\vbox{%
                 \noindent\hbox to \wd0{\hfill\GRAPHIC{#5}{#4}{#1}{#2}{#3}\hfill}\\%
                 \noindent\hbox{\QCBOptA}%
               }%
               \wd1=\wd0
            \fi
         \else
            \ifdim\wd0>0pt
              \hsize=\@tempdima
              \setbox1=\vbox{%
                \unskip\GRAPHIC{#5}{#4}{#1}{#2}{0pt}%
                \break
                \unskip\hbox to \@tempdima{\hfill \QCBOptA\hfill}%
              }%
              \wd1=\@tempdima
           \else
              \hsize=\@tempdima
              \setbox1=\vbox{%
                \unskip\GRAPHIC{#5}{#4}{#1}{#2}{0pt}%
              }%
              \wd1=\@tempdima
           \fi
         \fi
         \@tempdimb=\ht1
         \advance\@tempdimb by -#2
         \advance\@tempdimb by #3
         \leavevmode
         \raise -\@tempdimb \hbox{\box1}%
      \fi
      \egroup%
}%
\def\DFRAME#1#2#3#4#5{%
  \vspace\topsep
  \hfil\break
  \bgroup
     \leftskip\@flushglue
	 \rightskip\@flushglue
	 \parindent\z@
	 \parfillskip\z@skip
     \let\QCTOptA\empty
     \let\QCTOptB\empty
     \let\QCBOptA\empty
     \let\QCBOptB\empty
	 \vbox\bgroup
        \ifOverFrame
           #5\QCTOptA\par
        \fi
        \GRAPHIC{#4}{#3}{#1}{#2}{\z@}%
        \ifUnderFrame
           \break#5\QCBOptA
        \fi
	 \egroup
  \egroup
  \vspace\topsep
  \break
}%
\def\FFRAME#1#2#3#4#5#6#7{%
  \@ifundefined{floatstyle}
    {
     \begin{figure}[#1]%
    }
    {
	 \ifx#1h
      \begin{figure}[H]%
	 \else
      \begin{figure}[#1]%
	 \fi
	}
  \let\QCTOptA\empty
  \let\QCTOptB\empty
  \let\QCBOptA\empty
  \let\QCBOptB\empty
  \ifOverFrame
    #4
    \ifx\QCTOptA\empty
    \else
      \ifx\QCTOptB\empty
        \caption{\QCTOptA}%
      \else
        \caption[\QCTOptB]{\QCTOptA}%
      \fi
    \fi
    \ifUnderFrame\else
      \label{#5}%
    \fi
  \else
    \UnderFrametrue%
  \fi
  \begin{center}\GRAPHIC{#7}{#6}{#2}{#3}{\z@}\end{center}%
  \ifUnderFrame
    #4
    \ifx\QCBOptA\empty
      \caption{}%
    \else
      \ifx\QCBOptB\empty
        \caption{\QCBOptA}%
      \else
        \caption[\QCBOptB]{\QCBOptA}%
      \fi
    \fi
    \label{#5}%
  \fi
  \end{figure}%
 }%
\newcount\dispkind%

\def\makeactives{
  \catcode`\"=\active
  \catcode`\;=\active
  \catcode`\:=\active
  \catcode`\'=\active
  \catcode`\~=\active
}
\bgroup
   \makeactives
   \gdef\activesoff{%
      \def"{\string"}%
      \def;{\string;}%
      \def:{\string:}%
      \def'{\string'}%
      \def~{\string~}%
    }
\egroup

\def\FRAME#1#2#3#4#5#6#7#8{%
 \bgroup
 \ifnum\@msidraft=\@ne
   \wasdrafttrue
 \else
   \wasdraftfalse%
 \fi
 \def\LaTeXparams{}%
 \dispkind=\z@
 \def\LaTeXparams{}%
 \doFRAMEparams{#1}%
 \ifnum\dispkind=\z@\IFRAME{#2}{#3}{#4}{#7}{#8}{#5}\else
  \ifnum\dispkind=\@ne\DFRAME{#2}{#3}{#7}{#8}{#5}\else
   \ifnum\dispkind=\tw@
    \edef\@tempa{\noexpand\FFRAME{\LaTeXparams}}%
    \@tempa{#2}{#3}{#5}{#6}{#7}{#8}%
    \fi
   \fi
  \fi
  \ifwasdraft\@msidraft=1\else\@msidraft=0\fi{}%
  \egroup
 }%
\def\TEXUX#1{"texux"}

\def\limfunc#1{\mathop{\rm #1}}%
\def\func#1{\mathop{\rm #1}\nolimits}%
\long\def\QQQ#1#2{%
     \long\expandafter\def\csname#1\endcsname{#2}}%
\@ifundefined{QTP}{\def\QTP#1{}}{}
\@ifundefined{QEXCLUDE}{\def\QEXCLUDE#1{}}{}
\@ifundefined{Qlb}{}{}
\@ifundefined{Qlt}{}{}
\long\def\QQA#1#2{}%
\def\QTR#1#2{{\csname#1\endcsname {#2}}}%

\def\EXPAND#1[#2]#3{}%
\def\NOEXPAND#1[#2]#3{}%
\def\LaTeXparent#1{}%
\def\ChildStyles#1{}%
\def\ChildDefaults#1{}%
\def\QTagDef#1#2#3{}%

\@ifundefined{correctchoice}{}{}
\@ifundefined{HTML}{\def\HTML#1{\relax}}{}
\@ifundefined{TCIIcon}{\def\TCIIcon#1#2#3#4{\relax}}{}
\if@compatibility
\else
  \providecommand{\UNICODE}[2][]{\protect\rule{.1in}{.1in}}
  \providecommand{\U}[1]{\protect\rule{.1in}{.1in}}
  
\fi

\@ifundefined{lambdabar}{
      
   }{}

\@ifundefined{StyleEditBeginDoc}{}{}
\def\QQfnmark#1{\footnotemark}

\@ifundefined{TCIMAKEINDEX}{}{\makeindex}%
\@ifundefined{abstract}{%
 \def\abstract{%
  \if@twocolumn
   \section*{Abstract (Not appropriate in this style!)}%
   \else \small
   \begin{center}{\bf Abstract\vspace{-.5em}\vspace{\z@}}\end{center}%
   \quotation
   \fi
  }%
 }{%
 }%
\@ifundefined{endabstract}{\def\endabstract
  {\if@twocolumn\else\endquotation\fi}}{}%
\@ifundefined{maketitle}{\def\maketitle#1{}}{}%
\@ifundefined{affiliation}{\def\affiliation#1{}}{}%
\@ifundefined{proof}{}{}%
\@ifundefined{endproof}{}{}%
\@ifundefined{newfield}{\def\newfield#1#2{}}{}%
\@ifundefined{chapter}{\def\chapter#1{\par(Chapter head:)#1\par }%
 \newcount\c@chapter}{}%
\@ifundefined{part}{\def\part#1{\par(Part head:)#1\par }}{}%
\@ifundefined{section}{\def\section#1{\par(Section head:)#1\par }}{}%
\@ifundefined{subsection}{\def\subsection#1%
 {\par(Subsection head:)#1\par }}{}%
\@ifundefined{subsubsection}{\def\subsubsection#1%
 {\par(Subsubsection head:)#1\par }}{}%
\@ifundefined{paragraph}{\def\paragraph#1%
 {\par(Subsubsubsection head:)#1\par }}{}%
\@ifundefined{subparagraph}{\def\subparagraph#1%
 {\par(Subsubsubsubsection head:)#1\par }}{}%
\@ifundefined{therefore}{}{}%
\@ifundefined{backepsilon}{}{}%
\@ifundefined{yen}{}{}%
\@ifundefined{registered}{%
   \def\registered{\relax\ifmmode{}\r@gistered
                    \else$\m@th\r@gistered$\fi}%
 \def\r@gistered{^{\ooalign
  {\hfil\raise.07ex\hbox{$\scriptstyle\rm\text{R}$}\hfil\crcr
  \mathhexbox20D}}}}{}%
\@ifundefined{Eth}{}{}%
\@ifundefined{eth}{}{}%
\@ifundefined{Thorn}{}{}%
\@ifundefined{thorn}{}{}%
\@ifundefined{degree}{}{}%
\newdimen\theight
\@ifundefined{Column}{\def\Column{%
 \vadjust{\setbox\z@=\hbox{\scriptsize\quad\quad tcol}%
  \theight=\ht\z@\advance\theight by \dp\z@\advance\theight by \lineskip
  \kern -\theight \vbox to \theight{%
   \rightline{\rlap{\box\z@}}%
   \vss
   }%
  }%
 }}{}%
\@ifundefined{qed}{\def\qed{%
 \ifhmode\unskip\nobreak\fi\ifmmode\ifinner\else\hskip5\p@\fi\fi
 \hbox{\hskip5\p@\vrule width4\p@ height6\p@ depth1.5\p@\hskip\p@}%
 }}{}%
\@ifundefined{cents}{}{}%
\@ifundefined{tciLaplace}{}{}%
\@ifundefined{tciFourier}{}{}%
\@ifundefined{textcurrency}{}{}%
\@ifundefined{texteuro}{}{}%
\@ifundefined{euro}{}{}%
\@ifundefined{textfranc}{}{}%
\@ifundefined{textlira}{}{}%
\@ifundefined{textpeseta}{}{}%
\@ifundefined{miss}{\def\miss{\hbox{\vrule height2\p@ width 2\p@ depth\z@}}}{}%
\@ifundefined{vvert}{}{}
\@ifundefined{tcol}{\def\tcol#1{{\baselineskip=6\p@ \vcenter{#1}} \Column}}{}%
\@ifundefined{dB}{}{}
\@ifundefined{mB}{}{}
\@ifundefined{nB}{}{}
\@ifundefined{note}{}{}%
\def\newfmtname{LaTeX2e}
\ifx\fmtname\newfmtname
  \DeclareOldFontCommand{\rm}{\normalfont\rmfamily}{\mathrm}
  \DeclareOldFontCommand{\sf}{\normalfont\sffamily}{\mathsf}
  \DeclareOldFontCommand{\tt}{\normalfont\ttfamily}{\mathtt}
  \DeclareOldFontCommand{\bf}{\normalfont\bfseries}{\mathbf}
  \DeclareOldFontCommand{\it}{\normalfont\itshape}{\mathit}
  \DeclareOldFontCommand{\sl}{\normalfont\slshape}{\@nomath\sl}
  \DeclareOldFontCommand{\sc}{\normalfont\scshape}{\@nomath\sc}
\fi

\def\alpha{{\Greekmath 010B}}%
\def\beta{{\Greekmath 010C}}%
\def\gamma{{\Greekmath 010D}}%
\def\delta{{\Greekmath 010E}}%
\def\epsilon{{\Greekmath 010F}}%
\def\zeta{{\Greekmath 0110}}%
\def\eta{{\Greekmath 0111}}%
\def\theta{{\Greekmath 0112}}%
\def\iota{{\Greekmath 0113}}%
\def\kappa{{\Greekmath 0114}}%
\def\lambda{{\Greekmath 0115}}%
\def\mu{{\Greekmath 0116}}%
\def\nu{{\Greekmath 0117}}%
\def\xi{{\Greekmath 0118}}%
\def\pi{{\Greekmath 0119}}%
\def\rho{{\Greekmath 011A}}%
\def\sigma{{\Greekmath 011B}}%
\def\tau{{\Greekmath 011C}}%
\def\upsilon{{\Greekmath 011D}}%
\def\phi{{\Greekmath 011E}}%
\def\chi{{\Greekmath 011F}}%
\def\psi{{\Greekmath 0120}}%
\def\omega{{\Greekmath 0121}}%
\def\varepsilon{{\Greekmath 0122}}%
\def\vartheta{{\Greekmath 0123}}%
\def\varpi{{\Greekmath 0124}}%
\def\varrho{{\Greekmath 0125}}%
\def\varsigma{{\Greekmath 0126}}%
\def\varphi{{\Greekmath 0127}}%

\def\nabla{{\Greekmath 0272}}
\def\FindBoldGroup{%
   {\setbox0=\hbox{$\mathbf{x\global\edef\theboldgroup{\the\mathgroup}}$}}%
}

\def\Greekmath#1#2#3#4{%
    \if@compatibility
        \ifnum\mathgroup=\symbold
           \mathchoice{\mbox{\boldmath$\displaystyle\mathchar"#1#2#3#4$}}%
                      {\mbox{\boldmath$\textstyle\mathchar"#1#2#3#4$}}%
                      {\mbox{\boldmath$\scriptstyle\mathchar"#1#2#3#4$}}%
                      {\mbox{\boldmath$\scriptscriptstyle\mathchar"#1#2#3#4$}}%
        \else
           \mathchar"#1#2#3#4%
        \fi
    \else
        \FindBoldGroup
        \ifnum\mathgroup=\theboldgroup 
           \mathchoice{\mbox{\boldmath$\displaystyle\mathchar"#1#2#3#4$}}%
                      {\mbox{\boldmath$\textstyle\mathchar"#1#2#3#4$}}%
                      {\mbox{\boldmath$\scriptstyle\mathchar"#1#2#3#4$}}%
                      {\mbox{\boldmath$\scriptscriptstyle\mathchar"#1#2#3#4$}}%
        \else
           \mathchar"#1#2#3#4%
        \fi     	
	  \fi}

\newif\ifGreekBold  \GreekBoldfalse
\let\SAVEPBF=\pbf
\def\pbf{\GreekBoldtrue\SAVEPBF}%

\@ifundefined{theorem}{\newtheorem{theorem}{Theorem}}{}
\@ifundefined{lemma}{\newtheorem{lemma}[theorem]{Lemma}}{}
\@ifundefined{corollary}{\newtheorem{corollary}[theorem]{Corollary}}{}
\@ifundefined{conjecture}{}{}
\@ifundefined{proposition}{\newtheorem{proposition}[theorem]{Proposition}}{}
\@ifundefined{axiom}{}{}
\@ifundefined{remark}{\newtheorem{remark}{Remark}}{}
\@ifundefined{example}{\newtheorem{example}{Example}}{}
\@ifundefined{exercise}{}{}
\@ifundefined{definition}{\newtheorem{definition}{Definition}}{}

\@ifundefined{mathletters}{%
  \newcounter{equationnumber}
  \def\mathletters{%
     \addtocounter{equation}{1}
     \edef\@currentlabel{\theequation}%
     \setcounter{equationnumber}{\c@equation}
     \setcounter{equation}{0}%
     \edef\theequation{\@currentlabel\noexpand\alph{equation}}%
  }
  
}{}

\@ifundefined{BibTeX}{%
    \def\BibTeX{{\rm B\kern-.05em{\sc i\kern-.025em b}\kern-.08em
                 T\kern-.1667em\lower.7ex\hbox{E}\kern-.125emX}}}{}%
\@ifundefined{AmS}%
    {\def\AmS{{\protect\usefont{OMS}{cmsy}{m}{n}%
                A\kern-.1667em\lower.5ex\hbox{M}\kern-.125emS}}}{}%
\@ifundefined{AmSTeX}{}{}%
\def\@@eqncr{\let\@tempa\relax
    \ifcase\@eqcnt \def\@tempa{& & &}\or \def\@tempa{& &}%
      \else \def\@tempa{&}\fi
     \@tempa
     \if@eqnsw
        \iftag@
           \@taggnum
        \else
           \@eqnnum\stepcounter{equation}%
        \fi
     \fi
     \global\tag@false
     \global\@eqnswtrue
     \global\@eqcnt\z@\cr}

\def\TCItag{\@ifnextchar*{\@TCItagstar}{\@TCItag}}
\def\@TCItag#1{%
    \global\tag@true
    \global\def\@taggnum{(#1)}%
    \global\def\@currentlabel{#1}}
\def\@TCItagstar*#1{%
    \global\tag@true
    \global\def\@taggnum{#1}%
    \global\def\@currentlabel{#1}}
\def\tint{\msi@int\textstyle\int}%
\def\tiint{\msi@int\textstyle\iint}%
\def\tiiint{\msi@int\textstyle\iiint}%
\def\tiiiint{\msi@int\textstyle\iiiint}%
\def\tidotsint{\msi@int\textstyle\idotsint}%
\def\toint{\msi@int\textstyle\oint}%

\def\tsum{\mathop{\textstyle \sum }}%
\def\tbigcup{\mathop{\textstyle \bigcup }}%
\newtoks\temptoksa
\newtoks\temptoksb
\newtoks\temptoksc

\def\msi@int#1#2{%
 \def\@temp{{#1#2\the\temptoksc_{\the\temptoksa}^{\the\temptoksb}}}%
 \futurelet\@nextcs
 \@int
}

\def\@int{%
   \ifx\@nextcs\limits
      \typeout{Found limits}%
      \temptoksc={\limits}%
	  \let\@next\@intgobble%
   \else\ifx\@nextcs\nolimits
      \typeout{Found nolimits}%
      \temptoksc={\nolimits}%
	  \let\@next\@intgobble%
   \else
      \typeout{Did not find limits or no limits}%
      \temptoksc={}%
      \let\@next\msi@limits%
   \fi\fi
   \@next
}%

\def\@intgobble#1{%
   \typeout{arg is #1}%
   \msi@limits
}

\def\msi@limits{%
   \temptoksa={}%
   \temptoksb={}%
   \@ifnextchar_{\@limitsa}{\@limitsb}%
}

\def\@limitsa_#1{%
   \temptoksa={#1}%
   \@ifnextchar^{\@limitsc}{\@temp}%
}

\def\@limitsb{%
   \@ifnextchar^{\@limitsc}{\@temp}%
}

\def\@limitsc^#1{%
   \temptoksb={#1}%
   \@ifnextchar_{\@limitsd}{\@temp}%
}

\def\@limitsd_#1{%
   \temptoksa={#1}%
   \@temp
}

\def\dint{\msi@int\displaystyle\int}%
\def\diint{\msi@int\displaystyle\iint}%
\def\diiint{\msi@int\displaystyle\iiint}%
\def\diiiint{\msi@int\displaystyle\iiiint}%
\def\didotsint{\msi@int\displaystyle\idotsint}%
\def\doint{\msi@int\displaystyle\oint}%

\if@compatibility\else
  \RequirePackage{amsmath}
\fi

\begin{document}

\title{Relaxed Lagrangian duality in convex infinite optimization:
reducibility and strong duality}
\author{N. Dinh\thanks{%
International University, Vietnam National University - HCMC,
(ndinh@hcmiu.edu.vn)},\thanks{%
Vietnam National University - HCM city, Linh Trung ward, Thu Duc district,
Ho Chi Minh city, Vietnam }, \ \ M. A. Goberna\thanks{%
Department of Mathematics, University of Alicante, Alicante, Spain
(mgoberna@ua.es)}, \ \ M. A. L\'{o}pez\thanks{%
Department of Mathematics, University of Alicante, Alicante, Spain
(marco.antonio@ua.es); and CIAO, Federation University, Ballarat, Australia,
corresponding author}, \ \ M. Volle\thanks{%
Avignon University, LMA EA 2151, Avignon, France
(michel.volle@univ-avignon.fr)} }
\maketitle
\date{}

\begin{abstract}
We associate with each convex optimization problem, posed on some locally
convex space, {with infinitely many constraints indexed by the set} $T,$\
and a given non-empty family $\mathcal{H}$ {of} finite subsets of $T,$ a
suitable Lagrangian-Haar dual problem. We obtain necessary and sufficient
conditions for $\mathcal{H}$-reducibility, that is, equivalence to some
subproblem obtained by replacing the whole index set $T$\ by some element of
$\mathcal{H}$. Special attention is addressed to linear optimization,
infinite and semi-infinite, and to convex problems {with a countable family
of constraints}. Results on zero $\mathcal{H}$-duality gap and on $\mathcal{H%
}$-(stable) strong duality are provided. Examples are given along the paper
to illustrate the meaning of the results.
\end{abstract}




\section{Introduction}

The aim of this paper is to analyze relaxed forms of the \textit{%
Lagrangian-Haar dual problem} relative to the \textit{convex infinite
programming }(CIP, in brief) \textit{problem}
\begin{equation}
(\mathrm{P})\ \ \ \ \ \ \inf f(x)\ \ \text{ s.t.}\ \ f_{t}(x)\leq 0,\ \ t\in
T,  \label{1.1}
\end{equation}%
where $X$ is a locally convex Hausdorff topological vector space, $T$ is an
arbitrary infinite index set, and $\{f;\ f_{t},t\in T\}$ are convex proper
functions on $X$. The key stone in this paper is the\textbf{\ }(\textit{%
relaxed}) $\mathcal{H}$\textit{-Lagrangian-Haar dual} of $(\mathrm{P}),$
where $\mathcal{H}$ is a given family of nonempty finite subsets of the
index set $T,$ defined as
\begin{equation*}
(\mathrm{D}_{\mathcal{H}})\ \ \ \ \ \ \sup\limits_{H\in \mathcal{H},\ \mu
\in \mathbb{R}_{+}^{H}}\inf\limits_{x\in X}\left\{ f(x)+\sum\limits_{t\in
H}\mu _{t}f_{t}(x)\right\} ,
\end{equation*}%
where $\mu =(\mu _{t})_{t\in H}\in \mathbb{R}_{+}^{H},$ with the rule $%
0\times (+\infty )=0$, which {is applied }along the whole paper, except in
Remark \ref{rem2}. We say that \textit{zero }$\mathcal{H}$-\textit{duality}
(or just $\mathcal{H}$-\textit{duality}) holds if the optimal values of $(%
\mathrm{P})$ and $(\mathrm{D}_{\mathcal{H}})$ coincide, i.e., if $\inf (%
\mathrm{P})=\sup (\mathrm{D}_{\mathcal{H}}),$ while $\mathcal{H}$-\textit{%
strong duality} holds if, additionally, $(\mathrm{D}_{\mathcal{H}})$ is
solvable, i.e., if $\inf (\mathrm{P})=\max (\mathrm{D}_{\mathcal{H}}).$

The usual \textit{Lagrangian-Haar} dual $(\mathrm{D})$ of $(\mathrm{P})$
(see, e.g., \cite{DGLS07}, \cite{FLN09}, \cite{Volle2014}, \cite{Volle2015})
corresponds to the case where $\mathcal{H}$ is the family $\mathcal{F}(T)$
of all non-empty finite subsets of $T$, that is,
\begin{equation}
(\mathrm{D})\ \ \ \ \ \ \sup\limits_{H\in \mathcal{F}(T),\ \mu \in \mathbb{R}%
_{+}^{H}}\inf\limits_{x\in X}\left\{ f(x)+\sum\limits_{t\in H}\mu
_{t}f_{t}(x)\right\} .  \label{1.2}
\end{equation}%
{Obviously }$\sup (\mathrm{D}_{\mathcal{H}})\leq \sup (\mathrm{D}).$ {%
Moreover, the so-called weak duality inequality establishes that }$\sup (%
\mathrm{D})\leq \inf (\mathrm{P}).$

Other examples of such type of families are $\mathcal{H}_{1}:=\left\{
\{t\},\ t\in T\right\} $ and, for $T=\mathbb{N}$, $\mathcal{H}_{\mathbb{N}%
}:=\left\{ \{1,\ldots ,n\},\ n\in \mathbb{N}\right\} ,$ which are also {%
meaningful }in the framework of duality for the robust sum of functions \cite%
{DGV20}. So, we also pay attention to the dual problems
\begin{equation*}
(\mathrm{D}_{\mathcal{H}_{1}})\ \ \ \ \ \ \sup\limits_{(t,\mu )\in T\times
\mathbb{R}_{+}}\inf_{x\in X}\left\{ f(x)+\mu f_{t}(x)\right\} ,
\end{equation*}%
and%
\begin{equation*}
(\mathrm{D}_{\mathcal{H}_{\mathbb{N}}})\ \ \ \ \ \ \sup\limits_{n\in \mathbb{%
N}\text{, }\mu \in \mathbb{R}_{+}^{n}}\inf_{x\in X}\left\{
f(x)+\sum_{k=1}^{n}\mu _{k}f_{k}(x)\right\} .
\end{equation*}%
%
%
%
%
%
%
%
%
%
%
%
%
%
%
%
%
%
%
%
%
%
%
%
%
%
%
%
%
%
%
%

The problem $(\mathrm{P})$ in (\ref{1.1})\ is said to be \textit{reducible}
if there exists a finite set $H\subset T$ such that the optimal value of $(%
\mathrm{P})$ coincides with that of the problem $(\mathrm{P}_{H})$ which
results of replacing $T$ by $H$ in $(\mathrm{P}).$ The reducible linear
semi-infinite programming problems have been characterized in \cite[Theorem
8.3]{GL98}, but we do not know antecedents on reducibility for other types
of convex optimization problems. In the same vein, we can say that a dual
problem $(\mathrm{D})$ is \textit{reducible} if there exists a finite set $%
H\subset T$ such that the optimal value of $(\mathrm{D})$ coincides with
that of%
\begin{equation*}
(\mathrm{D}_{H})\ \ \ \ \ \ \sup\limits_{\mu \in \mathbb{R}%
_{+}^{H}}\inf\limits_{x\in X}\left\{ f(x)+\sum\limits_{t\in H}\mu
_{t}f_{t}(x)\right\} .
\end{equation*}%
Accordingly, we say that $(\mathrm{P})$ (respectively, $(\mathrm{D})$) is $%
\mathcal{H}$-\textit{reducible} if there exists a finite set $H\in \mathcal{H%
}$ such that its optimal value coincides with that of $(\mathrm{P}_{H})$ ({%
respectively}, $(\mathrm{D}_{H})$). So, a given CIP problem $(\mathrm{P})$
is $\mathcal{H}_{1}$-reducible if and only if it has the same optimal value
as some subproblem with a unique constraint, and it is $\mathcal{F}(T)$%
-reducible if and only if it is reducible (or $\mathcal{H}_{\mathbb{N}}$%
-reducible whenever $T$ is countable). Proposition 5.105 in \cite{BS00} can
be interpreted as providing a sufficient condition for a given convex
semi-infinite programming problem $(\mathrm{P}),$ with compact index set $T,$%
\ and its\ dual $(\mathrm{D})$, to be both reducible {relatively} to the
family $\mathcal{H}$ of subsets of $T$ whose cardinality is the dimension of
$X.$

While this paper is focused on $\mathcal{H}$-reducibility, zero $\mathcal{H}$%
-duality, $\mathcal{H}$-strong duality, and $\mathcal{H}$-stable strong
duality, in a forthcoming paper we consider reverse $\mathcal{H}$-strong
duality (where the solvable problem is $(\mathrm{P})$) and applications of
this type of relaxed Lagrangian duality to derive Farkas-type lemmas and
optimality conditions involving a fixed number of positive KKT multipliers
(e.g., one, if we choose $\mathcal{H}=\mathcal{H}_{1}$).

The classical Lagrange duality theory has been recently extended in another
direction in \cite{BBKY21}, replacing the convex functions by the so-called $%
\mathcal{H}$-convex functions, which are the supremum of certain class of
the space $\mathcal{H}$ of abstract affine functions (we use the same symbol
$\mathcal{H}$ for a class of index sets, so that our duality theorems are
independent of those of \cite{BBKY21}).

The paper is organized as follows. Section 2 reviews the classical duality
theorems for CIP problems and their finite subproblems. Section 3
characterizes the $\mathcal{H}$-reducibility of $(\mathrm{P})$ and $(\mathrm{%
D})$ in Theorem \ref{thm_reducible} in terms of $\mathcal{H}$-strong duality
of this pair of problems. In Section 4 special attention is addressed to $%
\mathcal{H}$-reducibility of linear infinite and semi-infinite programming
problems (in short LIP and LSIP, respectively). Section 5 is devoted to the $%
\mathcal{H}$-(stable) strong duality of the pair \textrm{(P)-}$(\mathrm{D}_{%
\mathcal{H}})$, which is characterized by Theorem \ref{thm2}. The very
particular case when $\mathcal{H}=\mathcal{H}_{1}$ is analyzed in Theorem %
\ref{thm1}. Also, a special attention is addressed to $\mathcal{H}_{\mathbb{N%
}}$-strong duality concerning countable convex infinite problems (Theorem %
\ref{thm3}, Corollary \ref{cor4}). Finally, Section 6 provides a
characterization of zero $\mathcal{H}$-duality gap between the problems
\textrm{(P)} and $(\mathrm{D}_{\mathcal{H}})$ (Theorem \ref{thm41a},
Corollaries \ref{corol61}-\ref{corol42b}).

\section{Preliminaries}

Let $X$ be a locally convex Hausdorff topological vector space, and suppose
that its\textbf{\ }topological dual $X^{\ast }$, with null element $%
0_{X^{\ast }},$ is endowed with the weak*-topology. The $w^{\ast }$-closure
of a set $A\subset X^{\ast }$ is denoted by $\overline{A}.$ If $\mathbb{A}%
\subset X^{\ast }\times \overline{\mathbb{R}},$ then $\overline{\mathbb{A}}$
denotes the closure of $\mathbb{A}$ w.r.t. the product topology. A subset $%
\mathbb{A}\subset X^{\ast }\times \mathbb{R}$ is said to be $w^{\ast }$%
-closed (respectively, $w^{\ast }$-closed convex) regarding another subset $%
\mathbb{B}\subset X^{\ast }\times \mathbb{R}$\ \ if \ \ $\bar{\mathbb{A}}%
\cap \mathbb{B}=\mathbb{A}\cap \mathbb{B}$ (respectively, $\left( \overline{%
\limfunc{co}}\mathbb{A}\right) \cap \mathbb{B}=\mathbb{A}\cap \mathbb{B}$),
see \cite{B10} (respectively, \cite{EV16}).

We denote by $\limfunc{co}A$ the convex hull of $A\subset X.$\ For a set $%
\emptyset \neq A\subset X$, by the convex\textbf{\ }cone generated by $A$ we
mean $\limfunc{cone}(A):=\mathbb{R}_{+}(\limfunc{co}A)=\{\mu x:\mu \in
\mathbb{R}_{+},\ x\in \limfunc{co}A\}.$

A function $h:X\rightarrow \overline{\mathbb{R}}:=\mathbb{R}\cup \{\pm
\infty \}$ is proper if its epigraph $\limfunc{epi}h$ is non-empty and never
takes the value $-\infty $; it is convex if $\limfunc{epi}h$ is convex; it
is lower semicontinuous (lsc, in brief) if $\limfunc{epi}h$ is closed; and
it is upper semicontinuous (usc, in brief) if $-h$ is lsc. For a proper
function $h,$ we denote by $[h\leq 0]:=\{x\in X:h(x)\leq 0\}$ its lower
level set of $0$ (and, similarly, the strict lower level set $[h<0]$)$,$ by $%
\func{dom}h\mathrm{,}$ $\limfunc{epi}\nolimits_{s}h,$ $\overline{h},$\ and $%
h^{\ast }$ its domain, its strict epigraph, its lsc envelope,\ and its
Legendre-Fenchel conjugate, respectively. We also denote by $\Gamma \left(
X\right) $ the class of lsc proper convex functions on $X$. By $\delta _{A}$
we denote the indicator function of $A\subset X,$ and $\delta _{A}\in \Gamma
\left( X\right) $ if and only if $A$ is closed, convex, and non-empty.

\subsection{Classical convex infinite optimization duality}

Now let $T$ be an infinite index set and $\{f;\ f_{t},t\in T\}$ be a
collection of convex proper functions. Consider the problem
\begin{equation*}
(\mathrm{P})\quad \inf f(x)\;\text{s.t.}\;f_{t}(x)\leq 0,\;t\in T,
\end{equation*}%
and its \textit{Lagrangian-Haar} dual, that is (\ref{1.2}) or, equivalently,
\begin{equation*}
(\mathrm{D})\quad \sup_{\lambda \in \mathbb{R}_{+}^{(T)}}\inf_{x\in
X}\left\{ f(x)+\left( \sum_{t\in T}\lambda _{t}f_{t}\right) (x)\right\} ,
\end{equation*}%
where%
\begin{equation*}
\mathbb{R}_{+}^{(T)}:=\{\lambda :T\rightarrow \mathbb{R}_{+}\text{ such that
}\lambda (t)\equiv \lambda _{t}=0\text{ for all }t\in T\text{ except for
finitely many}\}.
\end{equation*}%
Given\textbf{\ }$\lambda \in \mathbb{R}_{+}^{(T)}$, let us consider its
support,\textbf{\ } 
$\limfunc{supp}\lambda :=\{t\in T:\text{ }\lambda _{t}>0\}$,
and the associated function\textbf{\ }$\sum_{t\in T}\lambda _{t}f_{t}\colon
\ X\rightarrow \mathbb{R}\cup \{+\infty \},$ defined by\medskip
\begin{equation*}
\left( \sum_{t\in T}\lambda _{t}f_{t}\right) (x)=\left\{
\begin{tabular}{ll}
$\sum\limits_{t\in \limfunc{supp}\lambda }\lambda _{t}f_{t}(x),\text{ }$ & $%
\text{if }\limfunc{supp}\lambda \neq \emptyset ,$ \\
$0,$ & $\text{if }\limfunc{supp}\lambda =\emptyset \;\text{(i.e., $\lambda
=0_{T}\in \mathbb{R}_{+}^{(T)}$).}$%
\end{tabular}%
\right.
\end{equation*}%
\medskip

The following function\textbf{\ }$\varphi :X^{\ast }\rightarrow \overline{%
\mathbb{R}}$ is a key tool in our approach:\textbf{\ }%
\begin{equation*}
\varphi (x^{\ast }):=\inf_{\lambda \in \mathbb{R}_{+}^{(T)}}\left(
f+\sum_{t\in T}\lambda _{t}f_{t}\right) ^{\ast }(x^{\ast }).
\end{equation*}

Given the set in\textbf{\ }$X^{\ast }\times \mathbb{R}$\textbf{\ }%
\begin{equation*}
\mathcal{A}:=\bigcup_{\lambda \in \mathbb{R}_{+}^{(T)}}\limfunc{epi}\left(
f+\sum_{t\in T}\lambda _{t}f_{t}\right) ^{\ast },
\end{equation*}%
the following properties have been proved in \cite[(2.1) and (2.2)]%
{Volle2015} for $\{f;\ f_{t},t\in T\}\subset \Gamma \left( X\right) $, {\
although} they remain valid for arbitrary proper functions (even
non-convex):
\begin{equation}
\mathcal{A}\text{ is convex,\quad }\varphi \text{ is convex,\quad }\limfunc{%
epi}\nolimits_{s}\varphi \subset \mathcal{A}\subset \limfunc{epi}\varphi
\text{,\quad and }\limfunc{epi}\overline{\varphi }=\overline{\mathcal{A}}.
\label{marco7}
\end{equation}

Let
\begin{equation*}
E:=\bigcap_{t\in T}[f_{t}\leq 0]
\end{equation*}%
be the feasible set of $(\mathrm{P})$. Then, the weak duality relations
below always hold:
\begin{equation}
-\infty \leq (f+\delta _{E})^{\ast }(x^{\ast })\leq \varphi (x^{\ast })\leq
f^{\ast }(x^{\ast })\leq +\infty ,\;\forall x^{\ast }\in X^{\ast }.
\label{2.1bis}
\end{equation}%
In the case where $x^{\ast }=0_{X^{\ast }}$, one gets {from (\ref{2.1bis})}
\begin{equation*}
-\infty \leq \inf\nolimits_{X}f\leq \sup (\mathrm{D})\leq \inf (\mathrm{P}%
)\leq +\infty .
\end{equation*}

\subsection{Subprogram duality}

Given $H\in \mathcal{F}(T)$, consider the subproblem of $(\mathrm{P})$%
\medskip\
\begin{equation*}
(\mathrm{P}_{H})\quad \inf f(x)\;\text{s.t.}\;f_{t}(x)\leq 0,\;t\in H,
\end{equation*}%
\medskip and its Lagrangian dual
\begin{equation*}
(\mathrm{D}_{H})\quad \sup_{\mu \in \mathbb{R}_{+}^{H}}\inf_{x\in X}\left\{
f(x)+\left( \sum_{t\in H}\mu _{t}f_{t}\right) (x)\right\} ,
\end{equation*}%
\medskip where, for each $\mu \in \mathbb{R}_{+}^{H}$ and each $x\in X$,
\begin{equation*}
\left( \sum_{t\in H}\mu _{t}f_{t}\right) (x)=\sum_{t\in H}\mu _{t}f_{t}(x).
\end{equation*}

Define
\begin{equation}  \label{2.2b}
\mathcal{A}_{H}:=\bigcup\limits_{\mu \in \mathbb{R}_{+}^{H}}\limfunc{epi}%
\left( f+\sum_{t\in H}\mu _{t}f_{t}\right) ^{\ast }\subset X^{\ast }\times
\mathbb{R},
\end{equation}%
and, for each $x^{\ast }\in X^{\ast }$,
\begin{equation}  \label{2.2c}
\varphi _{H}(x^{\ast }):=\inf_{\mu \in \mathbb{R}_{+}^{H}}\left(
f+\sum_{t\in H}\mu _{t}f_{t}\right) ^{\ast }(x^{\ast })\in \overline{\mathbb{%
R}}.
\end{equation}%
We have that $\mathcal{A}_{H}$ is convex, $\varphi _{H}$ is convex, $%
\mathcal{A}_{H}\subset \mathcal{A}$, $\varphi _{H}\geq \varphi $, and $%
\limfunc{epi}_{s}\varphi _{H}\subset \mathcal{A}_{H}\subset \limfunc{epi}%
\varphi _{H}.$ Moreover,\textbf{\ }$\sup (\mathrm{D}_{H})\leq \sup (\mathrm{D%
})$ and, if\textbf{\ } 
$E_{H}:=\bigcap_{t\in H}[f_{t}\leq 0] $ 
is\textbf{\ }the feasible set of $(\mathrm{P}_{H})$, the next weak duality
relations hold:\medskip\
\begin{equation*}
-\infty \leq \left(f + \delta_E\right)^\ast (x^*) \leq (f+\delta
_{E_{H}})^{\ast }(x^{\ast })\leq \varphi _{H}(x^{\ast })\leq f^{\ast
}(x^{\ast })\leq +\infty ,\text{ }\forall x^{\ast }\in X^{\ast }.
\end{equation*}

\medskip For $x^{\ast }=0_{X^{\ast }}$ these relations yield
\begin{equation}
-\infty \leq \inf\nolimits_{X}f\leq \sup (\mathrm{D}_{H})\leq \inf (\mathrm{P%
}_{H})\leq \inf (\mathrm{P})\leq +\infty .  \label{marco52}
\end{equation}

\section{$\mathcal{H}$-reducibility}

Let $\mathcal{H}$ be a non-empty family of non-empty finite subsets of $T$,
i.e., $\emptyset \neq \mathcal{H}\subset \mathcal{F}(T)$. Recall that $(%
\mathrm{P})$ (respectively, $(\mathrm{D})$) is $\mathcal{H}$-reducible if
there exists $H\in \mathcal{H}$ such that $\inf (\mathrm{P})=\inf (\mathrm{P}%
_{H})$ (resp. $\sup (\mathrm{D})=\sup (\mathrm{D}_{H})$).

\subsection{Illustrative examples}

We now illustrate this desirable property with two examples.

\begin{example}
\label{Exam4.1}Consider the problem
\begin{equation}
\begin{tabular}{lll}
$\left( \mathrm{P}\right) $ & $\inf\limits_{x\in \mathbb{R}^{2}}$ & $f\left(
x\right) =\left\langle x^{\ast },x\right\rangle $ \\
& s.t. & $-tx_{1}+(t-1)x_{2}+t-t^{2}\leq 0,\;t\in \left[ 0,1\right] ,$%
\end{tabular}
\label{4.3}
\end{equation}%
where $x^{\ast }\in \mathbb{R}_{+}^{2}$, and its corresponding family of
singletons $\mathcal{H}_{1}=\left\{ \left\{ t\right\} :t\in \left[ 0,1\right]
\right\} .$ The feasible set $E$ of $\left( \mathrm{P}\right) $ is the
Minkowski sum of the curve $C=\left\{ \left( s,1+s-2\sqrt{s}\right) :s\in %
\left[ 0,1\right] \right\} $ with the positive cone $\mathbb{R}_{+}^{2}$
(see Figure 1). In fact each constraint of (\ref{4.3}) defines a supporting
halfspace to $C$ and viceversa. So, the boundary $\limfunc{bd}E$ of $E$ is
formed by those points of $E$ with a unique active constraint. \textbf{\ }%

\begin{center}
\includegraphics[width=1.95in]{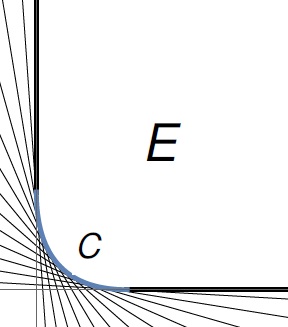}\\[1mm]
The feasible set of $\left(\mathrm{P}\right) $
\end{center}

If $x^{\ast }=\left( 0,0\right) ,$ given $H\in \mathcal{H}_{1},$\ we have
\begin{equation*}
\min \left( \mathrm{P}\right) =\min \left( \mathrm{P}_{H}\right) =0=\max
\left( \mathrm{D}_{H}\right) \leq \sup \left( \mathrm{D}\right) \leq \min
\left( \mathrm{P}\right) ,
\end{equation*}%
so that $(\mathrm{P})$ and $(\mathrm{D})$ are $\mathcal{H}_{1}$-reducible.
Hence, we can assume that $x^{\ast }\neq 0_{2}.$ Since any optimal solution $%
\overline{x}$, when it exists, lies on $\limfunc{bd}E,$ at such a point only
there will exist an active constraint, i.e., $\overline{t}\in \lbrack 0,1]$
such that
\begin{equation}
-\overline{t}\overline{x}_{1}+(\overline{t}-1)\overline{x}_{2}+\overline{t}-%
\overline{t}^{2}=0.  \label{4.5}
\end{equation}%
On the other hand, since the set $\left\{ (-t,t-1,t-t^{2}),\;t\in \left[ 0,1%
\right] \right\} $ is compact and any interior point of $E$ is a Slater
point, KKT optimality conditions (e.g. \cite[Theorems 7.1 and 5.3]{GL98})
apply to establish that $\overline{x}$ is optimal for problem (\ref{4.3}) if
and only if there exists $\lambda \geq 0$ such that
\begin{equation}
x^{\ast }=\lambda (\overline{t},1-\overline{t}).  \label{4.6}
\end{equation}%
So, the optimal solutions of $\left( \mathrm{P}\right) $ are those $x\in
\limfunc{bd}E$ such that there exist $\overline{t}\in \lbrack 0,1]$ and $%
\lambda \geq 0$ satisfying (\ref{4.5}) and (\ref{4.6}), which implies $%
x^{\ast }\in \mathbb{R}_{+}^{2}.$ \newline
If $x^{\ast }\notin \mathbb{R}_{+}^{2},$ one has $\inf \left( \mathrm{P}%
\right) =\sup (\mathrm{D})=-\infty $ and, recalling that the duality gap is
zero in linear programming (LP in short) except in the case when the primal
and the dual problems are simultaneously inconsistent, $\inf \left( \mathrm{P%
}_{H}\right) =\sup (\mathrm{D}_{H})=-\infty $ for all $H\in \mathcal{H}_{1}.$
So, $(\mathrm{P})$ and $(\mathrm{D})$ are $\mathcal{H}_{1}$-reducible. We
finally explore the non-trivial case that $x^{\ast }\in \mathbb{R}%
_{+}^{2}\diagdown \left\{ 0_{2}\right\} .$ Since
\begin{equation*}
\overline{x}=\left( \left( \frac{x_{2}^{\ast }}{x_{1}^{\ast }+x_{2}^{\ast }}%
\right) ^{2},\left( \frac{x_{1}^{\ast }}{x_{1}^{\ast }+x_{2}^{\ast }}\right)
^{2}\right) \in C,\text{ with }\overline{t}=\frac{x_{1}^{\ast }}{x_{1}^{\ast
}+x_{2}^{\ast }}\text{ and }\lambda =x_{1}^{\ast }+x_{2}^{\ast }>0,
\end{equation*}%
satisfies (\ref{4.5}) and (\ref{4.6}), $\overline{x}$ is an optimal solution
of $\left( \mathrm{P}\right) $ (the unique one when $x^{\ast }\in \mathbb{R}%
_{++}^{2},$ as Figure 1 shows), and $\min \left( \mathrm{P}\right) =\frac{%
x_{1}^{\ast }x_{2}^{\ast }}{x_{1}^{\ast }+x_{2}^{\ast }}.$ It is easy to see
that, taking $H:=\left\{ \overline{t}\right\} \in \mathcal{H}_{1},$ $%
\overline{x}$ is also optimal for $\left( \mathrm{P}_{H}\right) $ (but not
the unique one, as $\left( \mathrm{P}_{H}\right) $ has infinitely many
optimal solutions) and $\min \left( \mathrm{P}\right) =\min \left( \mathrm{P}%
_{H}\right) ;$ hence $\left( \mathrm{P}\right) $ is $\mathcal{H}_{1}$\textit{%
-reducible. It is also straightforward that }the unique feasible solution of
$\left( \mathrm{D}_{H}\right) $ is $\mu _{\overline{t}}=x_{1}^{\ast
}+x_{2}^{\ast }>0,$ with $\max \left( \mathrm{D}_{H}\right) =\frac{%
x_{1}^{\ast }x_{2}^{\ast }}{x_{1}^{\ast }+x_{2}^{\ast }}.$ Since $\max
\left( \mathrm{D}_{H}\right) \leq \sup \left( \mathrm{D}\right) \leq \min
\left( \mathrm{P}\right) =\min \left( \mathrm{P}_{H}\right) $ and $\max
\left( \mathrm{D}_{H}\right) =\min \left( \mathrm{P}_{H}\right) $, the
problem $\left( \mathrm{D}\right) $ is also $\mathcal{H}_{1}$\textit{%
-reducible. }

From the previous discussion we conclude that \textit{both problems, }$(%
\mathrm{P})$ and $(\mathrm{D}),$ $\mathcal{H}_{1}$\textit{-reducible}$.$
\end{example}

\begin{example}
\label{Exam4.2} Consider the LSIP problem
\begin{equation*}
\begin{tabular}{lll}
$\left( \mathrm{P}\right) $ & $\inf\limits_{x\in \mathbb{R}^{2}}$ & $f\left(
x\right) =\left\langle x^{\ast },x\right\rangle $ \\
& s.t. & $-tx_{1}+(t-1)x_{2}-t^{2}+t-\frac{1}{r}\leq 0,\;\left( t,r\right)
\in \left[ 0,1\right] \times \mathbb{N},$%
\end{tabular}%
\end{equation*}%
where $x^{\ast }\in \mathbb{R}_{+}^{2}$ and the family of singletons is now $%
\mathcal{H}_{1}=\left\{ \left\{ \left( t,r\right) \right\} :t\in \left[ 0,1%
\right] ,r\in \mathbb{N}\right\} .$ The primal feasible set is still the set
$E$ in Figure 1, but now no constraint is active at the boundary of $E.$
Hence, the optimal value of the subproblem with a unique constraint
\begin{equation*}
\begin{tabular}{ll}
$\inf\limits_{x\in \mathbb{R}^{2}}$ & $f\left( x\right) =\left\langle
x^{\ast },x\right\rangle $ \\
s.t. & $-tx_{1}+(t-1)x_{2}-t^{2}+t-\frac{1}{r}\leq 0,$%
\end{tabular}%
\end{equation*}%
is less than $\min \left( \mathrm{P}\right) $ for any $x^{\ast }\in \mathbb{R%
}_{+}^{2}$ and $\left\{ \left( t,r\right) \right\} \in \mathcal{H}_{1},$ so
that $\left( \mathrm{P}\right) $ is not $\mathcal{H}_{1}$-reducible.
\end{example}

\subsection{ $\mathcal{H}$-reducibility of general convex problems}

The next two sets will play crucial roles along this paper. The first one,%
\begin{equation}
\mathcal{A}_{\mathcal{H}}:=\bigcup\limits_{H\in \mathcal{H}}\mathcal{A}%
_{H}=\bigcup\limits_{H\in \mathcal{H}}\bigcup\limits_{\mu \in \mathbb{R}%
_{+}^{H}}\limfunc{epi}\left( f+\sum_{t\in H}\mu _{t}f_{t}\right) ^{\ast },
\label{2.2bb}
\end{equation}%
involves all the data of $(\mathrm{P})$ while the second one,
\begin{equation*}
\mathcal{K}_{\mathcal{H}}:=\bigcup\limits_{H\in \mathcal{H}}\limfunc{cone}%
\left( \bigcup\limits_{t\in H}\limfunc{epi}f_{t}^{\ast }\right) ,
\end{equation*}%
only involves the constraints. In particular, $\mathcal{K}:=\mathcal{K}_{%
\mathcal{F}(T)}$ is convex and satisfies \cite[Theorem 3.1]{DGL06}
\begin{equation*}
E\neq \emptyset \Longleftrightarrow \left( 0_{X^{\ast }},-1\right) \notin
\overline{\mathcal{K}}.
\end{equation*}

{We also associate with} $\mathcal{H}$ {the function }
\begin{equation}
\varphi _{\mathcal{H}}:=\inf\limits_{H\in \mathcal{H}}\varphi
_{H}=\inf\limits_{H\in \mathcal{H}}\inf_{\mu \in \mathbb{R}_{+}^{H}}\left(
f+\sum_{t\in H}\mu _{t}f_{t}\right) ^{\ast }  \label{3.3b}
\end{equation}%
{and the set }%
\begin{equation*}
E_{\mathcal{H}}:=\bigcap\limits_{H\in \mathcal{H}}E_{H}=\bigcap\limits_{H\in
\mathcal{H}}\bigcap_{t\in H}[f_{t}\leq 0].
\end{equation*}%
{Obviously\ }$A_{\mathcal{H}}\subset \mathcal{A},\;\varphi _{\mathcal{H}%
}\geq \varphi ,\;$ $\limfunc{epi}_{s}\varphi _{\mathcal{H}}\subset \mathcal{A%
}_{\mathcal{H}}\subset \limfunc{epi}\varphi _{\mathcal{H}}$ and
\begin{equation}
\left( f+\delta _{E}\right) ^{\ast }(x^{\ast })\leq \inf\limits_{H\in
\mathcal{H}}(f+\delta _{E_{H}})^{\ast }(x^{\ast })\leq \varphi _{\mathcal{H}%
}(x^{\ast })\leq f^{\ast }(x^{\ast }),\text{ }\forall x^{\ast }\in X^{\ast },
\label{3.3k}
\end{equation}%
entailing{\ }
\begin{equation}
\sup (\mathrm{D}_{\mathcal{H}})\leq \sup (\mathrm{D}_{\mathcal{F}(T)})\equiv
\sup (\mathrm{D})\leq \inf (\mathrm{P}).  \label{3.3j}
\end{equation}

\begin{remark}
Let $\mathcal{H},\mathcal{J}\subset \mathcal{F}(T)$. If for each $H\in
\mathcal{H}$ there exists $J\in \mathcal{J}$ such that $H\subset J$, then $%
\mathcal{A}_{\mathcal{H}}\subset \mathcal{A}_{\mathcal{J}}$, $\varphi _{%
\mathcal{J}}\leq \varphi _{\mathcal{H}}$, and $\sup (\mathrm{D}_{\mathcal{H}%
})\leq \sup (\mathrm{D}_{\mathcal{J}}).$ This is, in particular, the case
when $\mathcal{H}\subset \mathcal{J}$.
\end{remark}

Let us examine the relationship {among } $\varphi _{\mathcal{H}}$, $\mathcal{%
A}_{\mathcal{H}}$ and $\mathcal{K}_{\mathcal{H}}$. Assume that $\{f;\
f_{t},t\in T\}\subset \Gamma (X)$ and\textbf{\ }$E\cap (\func{dom}f)\neq
\emptyset $, entailing that all the functions $f+\tsum\nolimits_{t\in H}\mu
_{t}f_{t}$ are proper. Then, following a similar reasoning to that in \cite[%
Section 3]{Volle2015}, we can establish, that
\begin{equation*}
\limfunc{epi}\overline{\varphi _{\mathcal{H}}}=\overline{\limfunc{epi}%
f^{\ast }+\mathcal{K}_{\mathcal{H}}}=\overline{\mathcal{A}_{\mathcal{H}}}.
\end{equation*}%
In the particular case of $f\equiv 0$ we get $\overline{\mathcal{A}_{%
\mathcal{H}}}=\overline{(\left\{ 0_{X^{\ast }}\right\} \times \mathbb{R}%
_{+})+\mathcal{K}_{\mathcal{H}}}=\overline{\mathcal{K}_{\mathcal{H}}}.$

We now characterize the $\mathcal{H}$-reducibility of a given bounded CIP
problem $\left( \mathrm{P}\right) $ by combining three main ingredients: the
set $\mathcal{A}_{\mathcal{H}},$ the strong duality between $\left( \mathrm{P%
}\right) $ and $(\mathrm{D}_{\mathcal{H}}\mathrm{),}$ and the following
property:%
\begin{equation*}
\left( \mathrm{R}_{\mathcal{H}}\right) \ \text{Each bounded subprogram }(%
\mathrm{P_{H}}),\ H\in \mathcal{H}\text{,}\mathcal{\ }\text{satisfies }\inf (%
\mathrm{P}_{H})=\max (\mathrm{D}_{H}).
\end{equation*}

Property $\left( \mathrm{R}_{\mathcal{H}}\right) $ is, for instance,
satisfied in the two well-known situations below:

\begin{enumerate}
\item $(\mathrm{P})$ is a LSIP problem: due to the fundamental LP-duality
theorem which establishes that if one of the problems $(\mathrm{P}_{H})$ or $%
(\mathrm{D}_{H})$ is bounded, then both problems are solvable and their
optimal values coincide.

\item $(\mathrm{P})$ is a general CIP problem such that (see \cite[Theorem
2.9.3]{Z02}):%
\begin{equation*}
\forall H\in \mathcal{H}\text{, }\exists a_{H}\in \func{dom}f\text{ such
that }f_{t}\left( a_{H}\right) <0,\text{ }\forall t\in H.
\end{equation*}%
This is in particular satisfied under the (generalized) Slater condition;
i.e.,%
\begin{equation*}
\exists a\in \func{dom}f\text{ such that }f_{t}\left( a\right) <0\text{ for
all }t\in T.
\end{equation*}
\end{enumerate}

\begin{theorem}[CIP $\mathcal{H}$-reducibility]
\label{thm_reducible} Let $($\textrm{$P$}$)$ be bounded (i.e., $\alpha
:=\inf ($\textrm{$P$}$)\in \mathbb{R}$) and consider the following
assertions involving a family $\emptyset \neq \mathcal{H}\subset \mathcal{F}%
(T):$\newline
$\mathrm{(i)}$ $\left( 0_{X^{\ast }},-\alpha \right) \in \mathcal{A}_{%
\mathcal{H}}.$\newline
$\mathrm{(ii)}$ $\mathcal{H}$-strong duality holds (i.e., $\alpha =\max (%
\mathrm{D}_{\mathcal{H}})$). \newline
$\mathrm{(iii)}$ $(\mathrm{D})$ is $\mathcal{H}$-reducible and zero $%
\mathcal{H}$-duality holds (entailing $\alpha =\sup (\mathrm{D})$).\newline
$\mathrm{(iv)}$ $(\mathrm{P})$ is $\mathcal{H}$-reducible. \newline
Then, one has $\mathrm{(i)}\Longleftrightarrow \mathrm{(ii)}\Longrightarrow
\mathrm{(iii)}\Longrightarrow \mathrm{(iv)}.$ Moreover, if assumption $%
\left( \mathrm{R}_{\mathcal{H}}\right) $ is satisfied, then $\mathrm{(i)},%
\mathrm{(ii)},\mathrm{(iii)},$ and $\mathrm{(iv)}$ are equivalent.
\end{theorem}

\noindent \textbf{Proof }$\left[ \mathrm{(i)}\Longrightarrow \mathrm{(ii)}%
\right] $ By definition of $\mathcal{A}_{\mathcal{H}}$ there exists $H\in
\mathcal{H}$ and $\mu \in \mathbb{R}_{+}^{H}$ such that \newline
$\left( f+\sum_{t\in H}\mu _{t}f_{t}\right) ^{\ast }\left( 0_{X^{\ast
}}\right) \leq -\alpha $ and, by (\ref{marco52}), we have%
\begin{equation*}
\alpha =\inf (\mathrm{P})\geq \inf (\mathrm{P}_{H})\geq \sup (\mathrm{D}%
_{H})\geq \inf\limits_{X}\left( f+\sum_{t\in H}\mu _{t}f_{t}\right) \geq
\alpha ,
\end{equation*}%
and%
\begin{equation*}
\sup (\mathrm{D})\leq \inf (\mathrm{P})=\max (\mathrm{D}_{H})\leq \sup (%
\mathrm{D}),
\end{equation*}%
which shows\textbf{\ }that $\mathrm{(ii)}$ holds.

$\left[ \mathrm{(ii)}\Longrightarrow \mathrm{(i)}\right] $ It is obvious
from the definition of $\mathcal{A}_{\mathcal{H}}.$

$\left[ \mathrm{(ii)}\Longrightarrow \mathrm{(iii)}\right] $ There exists $%
H\in \mathcal{H}$ and $\mu \in \mathbb{R}_{+}^{H}$ such that
\begin{equation*}
\sup (\mathrm{D})\leq \alpha =\inf\limits_{X}\left( f+\sum_{t\in H}\mu
_{t}f_{t}\right) \leq \sup (\mathrm{D}_{H})\leq \sup (\mathrm{D}),
\end{equation*}%
and $\mathrm{(iii)}$ holds.

$\left[ \mathrm{(iii)}\Longrightarrow \mathrm{(iv)}\right] $ There exists $%
H\in \mathcal{H}$ such that
\begin{equation*}
\alpha =\sup (\mathrm{D})=\sup (\mathrm{D}_{H})\leq \inf (\mathrm{P}%
_{H})\leq \alpha ,
\end{equation*}%
and $\mathrm{(iv)}$ holds.

$\left[ \mathrm{(iv)\Longrightarrow (ii)}\right] $ (under\textbf{\ }$\left(
\mathrm{R}_{\mathcal{H}}\right) $) There exists $H\in \mathcal{H}$ such
that, by $\left( \mathrm{R}_{\mathcal{H}}\right) ,$
\begin{equation*}
\alpha =\inf (\mathrm{P}_{H})=\max (\mathrm{D}_{H})\leq \sup (\mathrm{D}_{%
\mathcal{H}})\leq \sup (\mathrm{D})\leq \inf (\mathrm{P})=\alpha ,
\end{equation*}%
and $\mathrm{(ii)}$ holds$\mathrm{.}$\hfill $\square \medskip $

We need some new notions in order to establish new sufficient conditions for
$\mathcal{H}$-reducibility.

\begin{definition}
$\left( \mathrm{i}\right) $ A family $\mathcal{H}\subset \mathcal{F}(T)$ is
said to be \textit{covering} if $\bigcup\nolimits_{H\in \mathcal{H}}H=T$.%
\newline
$\left( \mathrm{ii}\right) $ A family $\mathcal{H}\subset \mathcal{F}(T)$ is
said to be \textit{directed} if for each $H,K\in \mathcal{H}$ there exists $%
L\in \mathcal{H}$ such that $H\cup K\subset L$.
\end{definition}

The families $\mathcal{F}(T)$ and $\mathcal{H}_{\mathbb{N}}$ are both
covering and directed, whereas\textbf{\ }$\mathcal{H}_{1}$ is a covering
family but it is not a directed one. Note that
\begin{equation*}
\mathcal{H}\text{ covering }\Longrightarrow \left[ E_{\mathcal{H}}=E\text{
and }\sup (\mathrm{D}_{\mathcal{H}_{1}})\leq \sup (\mathrm{D}_{\mathcal{H}})%
\right] .
\end{equation*}

\begin{proposition}
\label{pro1a} The following statements are equivalent:\newline
$\mathrm{(i)}$ $\mathcal{H}$ is a directed covering family.\newline
$\mathrm{(ii)}$ For each $K\in \mathcal{F}(T),$ there exists $H\in \mathcal{H%
}$ such that $K\subset H$.
\end{proposition}

\noindent \textbf{Proof } $[\mathrm{(i)}\Longrightarrow \mathrm{(ii)}]$ Let $%
K\in \mathcal{F}(T)$. Since $\mathcal{H}$ is covering one has $K\subset
\bigcup\limits_{t\in K}H_{t}$ with $t\in H_{t}\in \mathcal{H}$ for each $%
t\in K$. Since $K$ is finite and $\mathcal{H}$ is directed there exists $%
H\in \mathcal{H}$ such that $\bigcup\nolimits_{t\in K}H_{t}\subset H$ and $%
K\subset H$.

$[\mathrm{(ii)}\Longrightarrow \mathrm{(i)}]$ Let $t\in T$. By $\mathrm{(ii)}
$, with $K=\{t\}$, there exists $H\in \mathcal{H}$ such that $t\in H$, so
that $\mathcal{H}$ is covering. We prove now that $\mathcal{H}$ is directed.
Let $I,J\in \mathcal{H}$. By $\mathrm{(ii)}$ with $K=I\cup J$, there exists $%
H\in \mathcal{H}$ such that $I\cup J\subset H$. \hfill $\square $

\begin{proposition}
\label{pro1} For each directed covering family $\mathcal{H}\subset \mathcal{F%
}(T)$ one has
\begin{equation*}
\mathcal{A}_{\mathcal{H}}=\mathcal{A}_{\mathcal{F}(T)}=\mathcal{A},
\end{equation*}%
and, consequently,
\begin{equation*}
\varphi _{\mathcal{H}}=\varphi _{\mathcal{F}(T)}=\varphi ,\;\text{\ and }%
\sup (\mathrm{D}_{\mathcal{H}})=\sup (\mathrm{D}_{\mathcal{F}(T)})\equiv
\sup \mathrm{(D}).
\end{equation*}%
Moreover, $\mathcal{K}_{\mathcal{H}}=\mathcal{K}.$
\end{proposition}

\noindent \textbf{Proof } Since $\mathcal{A}_{\mathcal{H}}\subset \mathcal{A}%
_{\mathcal{F}(T)}\subset \mathcal{A}$, it suffices to prove that $\mathcal{A}%
\subset \mathcal{A}_{\mathcal{H}}$. Let $(x^{\ast },r)\in \mathcal{A}$.
There exists $\lambda \in \mathbb{R}_{+}^{(T)}$ such that $\left(
f+\sum_{t\in T}\lambda _{t}f_{t}\right) ^{\ast }(x^{\ast })\leq r$. If $%
\limfunc{supp}\lambda =\emptyset $ then $(x^{\ast },r)\in \limfunc{epi}%
f^{\ast },$ with $\limfunc{epi}f^{\ast }\subset \mathcal{A}_{H}$ for all $%
H\in \mathcal{H}$, and so, $(x^{\ast },r)\in \mathcal{A}_{\mathcal{H}}$.
Assume now that $K:=\limfunc{supp}\lambda \neq \emptyset $. By Proposition %
\ref{pro1a} there exists $H\in \mathcal{H}$ such that $K\subset H$ and we
have
\begin{equation*}
(x^{\ast },r)\in \mathcal{A}_{K}\subset \mathcal{A}_{H}\subset \mathcal{A}_{%
\mathcal{H}}.
\end{equation*}%
The proof of the relation $\mathcal{K}_{\mathcal{H}}=\mathcal{K}$ is
similar. \hfill $\square $

\begin{proposition}
\label{pro1c} Assume that $\mathcal{H}$ is directed covering. Then $\mathrm{%
(P)}$ (resp. $\mathrm{(D)}$) is reducible if and only if $\mathrm{(P)}$
(resp. $\mathrm{(D)}$) is $\mathcal{H}$-reducible.
\end{proposition}

\noindent \textbf{Proof } It suffices to prove that the necessity holds. So,
assume that $\mathrm{(P)}$ (resp. $\mathrm{(D)}$) is reducible. Then, there
exists $K\in \mathcal{F}(T)$ such that $\inf \mathrm{(P)}=\inf (\mathrm{P}%
_{K})$ (resp. $\sup \mathrm{(D)}=\sup (\mathrm{D}_{K})$). By Proposition \ref%
{pro1a} there exists $H\in \mathcal{H}$ such that $K\subset H$.
Consequently,
\begin{equation*}
\inf \mathrm{(P)}\geq \inf (\mathrm{P}_{H})\geq \inf (\mathrm{P}_{K})=\inf
\mathrm{(P)},
\end{equation*}%
(resp. $\sup \mathrm{(D)}\geq \sup (\mathrm{D}_{H})\geq \sup (\mathrm{D}%
_{K})=\sup \mathrm{(D)}$), and finally, $\inf \mathrm{(P)}=\inf (\mathrm{P}%
_{H})$ (resp. $\sup \mathrm{(D)}=\sup (\mathrm{D}_{H})$). \hfill $\square
\medskip $

\section{$\mathcal{H}$-reducibility of linear optimization problems}


We now consider the case where
\begin{equation}
f(x)=\left\langle x^{\ast },x\right\rangle ,\ f_{t}(x)=\left\langle
a_{t}^{\ast },x\right\rangle -b_{t},\;t\in T,  \label{4.2}
\end{equation}%
with $\left\{ x^{\ast };a_{t}^{\ast },t\in T\right\} \subset X^{\ast },$ and
$\left\{ b_{t},t\in T\right\} \subset \mathbb{R}.$ Equivalently,
\begin{equation}
(\mathrm{P})\ \ \ \ \ \ \inf \ \left\langle x^{\ast },x\right\rangle \ \
\text{ s.t.}\ \left\langle a_{t}^{\ast },x\right\rangle \leq b_{t},\ \ t\in
T,  \label{LIP}
\end{equation}%
and
\begin{equation*}
(\mathrm{D})\ \ \ \ \ \sup -\sum\limits_{t\in T}\lambda _{t}b_{t}\text{ s.t.
}\lambda \in \mathbb{R}_{+}^{(T)},~\sum\limits_{t\in T}\lambda
_{t}a_{t}^{\ast }=-x^{\ast },
\end{equation*}%
and, for each $H\in \mathcal{H}\subset \mathcal{F}(T)$,%
\begin{equation*}
(\mathrm{P}_{H})\ \ \ \ \ \ \inf \ \left\langle x^{\ast },x\right\rangle \ \
\text{ s.t.}\ \left\langle a_{t}^{\ast },x\right\rangle \leq b_{t},\ \ t\in
H,
\end{equation*}%
and%
\begin{equation*}
(\mathrm{D}_{H})\ \ \ \ \ \sup -\sum\limits_{t\in H}\mu _{t}b_{t}\ \ \text{
s.t. }\mu \in \mathbb{R}_{+}^{H},~\sum\limits_{t\in H}\mu _{t}a_{t}^{\ast
}=-x^{\ast },
\end{equation*}%
while%
\begin{equation*}
(\mathrm{D}_{\mathcal{H}})\ \ \ \ \ \sup -\sum\limits_{t\in H}\mu _{t}b_{t}\
\ \text{ s.t. }H\in \mathcal{H}\text{,\ }\mu \in \mathbb{R}%
_{+}^{H},~\sum\limits_{t\in H}\mu _{t}a_{t}^{\ast }=-x^{\ast }.
\end{equation*}%
Let us make explicit $\mathcal{A}_{H}$. We have%
\begin{eqnarray*}
\left( f+\tsum\nolimits_{t\in H}\mu _{t}f_{t}\right) ^{\ast } &=&\left(
x^{\ast }+\tsum\nolimits_{t\in H}\mu _{t}(a_{t}^{\ast }-b_{t})\right) ^{\ast
} \\
&=&\delta _{\left\{ x^{\ast }+\tsum\nolimits_{t\in H}\mu _{t}a_{t}^{\ast
}\right\} }+\tsum\nolimits_{t\in H}\mu _{t}b_{t},
\end{eqnarray*}%
\begin{equation*}
\limfunc{epi}\left( f+\tsum\nolimits_{t\in H}\mu _{t}f_{t}\right) ^{\ast
}=\left\{ (x^{\ast },0)\right\} +\left\{ \tsum\nolimits_{t\in H}\mu
_{t}(a_{t}^{\ast },b_{t})\right\} +\left\{ 0_{X^{\ast }}\right\} \times
\mathbb{R}_{+},
\end{equation*}%
and%
\begin{equation*}
\mathcal{A}_{H}=\tbigcup\limits_{\mu \in \mathbb{R}_{+}^{H}}\limfunc{epi}%
\left( f+\tsum\nolimits_{t\in H}\mu _{t}f_{t}\right) ^{\ast }=\left\{
(x^{\ast },0)\right\} +\mathcal{K}_{H},
\end{equation*}%
where
\begin{equation*}
\mathcal{K}_{H}:=\limfunc{cone}\left( \left\{ (a_{t}^{\ast },b_{t}),\ t\in
H\right\} +\left\{ 0_{X^{\ast }}\right\} \times \mathbb{R}_{+}\right) .
\end{equation*}

Taking into account that
\begin{equation}
\mathcal{K}_{\mathcal{H}}=\tbigcup\limits_{H\in \mathcal{H}}\mathcal{K}_{H},
\label{KH}
\end{equation}%
which is a cone (not necessarily convex), we get
\begin{equation}
\mathcal{A}_{\mathcal{H}}=\left\{ \left( x^{\ast },0\right) \right\} +%
\mathcal{K}_{\mathcal{H}}.  \label{4.1bis}
\end{equation}

\begin{proposition}
\label{Prop4.1}If $\mathcal{H}$ is a directed covering family, then%
\begin{equation*}
\mathcal{K}_{\mathcal{H}}=\limfunc{cone}\left( \left\{ (a_{t}^{\ast
},b_{t}),\ t\in T\right\} +\left\{ 0_{X^{\ast }}\right\} \times \mathbb{R}%
_{+}\right) .
\end{equation*}
\end{proposition}

\noindent \textbf{Proof }From Proposition \ref{pro1} and (\ref{KH}),
\begin{eqnarray*}
\mathcal{K}_{\mathcal{H}} &=&\mathcal{K}=\tbigcup\limits_{H\in \mathcal{F}%
(T)}\limfunc{cone}\left( \left\{ (a_{t}^{\ast },b_{t}),\ t\in H\right\}
+\left\{ 0_{X^{\ast }}\right\} \times \mathbb{R}_{+}\right) \\
&=&\tbigcup\limits_{H\in \mathcal{F}(T)}\tbigcup\limits_{\lambda \in \mathbb{%
R}_{+}^{(T)}}\left\{ \tsum\limits_{t\in H}\lambda _{t}(a_{t}^{\ast
},b_{t})+\left\{ 0_{X^{\ast }}\right\} \times \mathbb{R}_{+}\right\} \\
&=&\limfunc{cone}\left( \left\{ (a_{t}^{\ast },b_{t}),\ t\in T\right\}
+\left\{ 0_{X^{\ast }}\right\} \times \mathbb{R}_{+}\right) .
\end{eqnarray*}%
\hfill \hfill\ \hfill \hfill So, we are done.\hfill $\square $

\begin{corollary}[LIP $\mathcal{H}$-reducibility]
\label{Cor.4.2} Assume that $\left( \mathrm{P}\right) $ is bounded (that is,
$\alpha :=\inf ($\textrm{$P$}$)\in \mathbb{R}$), and
\begin{equation}
\forall H\in \mathcal{H}\text{, }\exists a_{H}\in X\text{ s.t. }\left\langle
a_{t}^{\ast },a_{H}\right\rangle <b_{t},\ \forall t\in H.  \label{4.2bis}
\end{equation}%
Then the following statements are equivalent:\newline
$\mathrm{(i)}$ $-\left( x^{\ast },\alpha \right) \in \mathcal{K}_{\mathcal{H}%
}.$\newline
$\mathrm{(ii)}$ $\mathcal{H}$-strong duality holds (i.e., $\alpha =\max (%
\mathrm{D}_{\mathcal{H}})$).\newline
$\mathrm{(iii)}$ $(\mathrm{D})$ is $\mathcal{H}$-reducible and zero $%
\mathcal{H}$-duality holds (entailing $\alpha =\sup (\mathrm{D})$).\newline
$\mathrm{(iv)}$ $(\mathrm{P})$ is $\mathcal{H}$-reducible.\newline
\end{corollary}

\noindent \textbf{Proof }By (\ref{4.2bis}), property $\left( \mathrm{R}_{%
\mathcal{H}}\right) $ is satisfied and we have just to apply Theorem \ref%
{thm_reducible} and use (\ref{4.1bis}). \hfill $\square $\bigskip

We now consider the LSIP problem $\left( \mathrm{P}\right) $ in (\ref{LIP})
with $X=\mathbb{R}^{n}.$

\begin{corollary}[LSIP $\mathcal{H}$-reducibility]
\label{Cor.4.3} Assume that $X=\mathbb{R}^{n}$ and $\left( \mathrm{P}\right)
$ is bounded. Then, the statements from $\mathrm{(i)}$ to $\mathrm{(iv)}$ in
Corollary \ref{Cor.4.2} are still equivalent. Moreover, if $\mathcal{H}$ is
a directed covering family, then the following statement is also equivalent
to the four mentioned statements:\newline
$\mathrm{(v)}$ $\alpha =\max (\mathrm{D})$ (usual strong duality).
\end{corollary}

\noindent \textbf{Proof }Since every LSIP problem satisfies $\left( \mathrm{R%
}_{\mathcal{H}}\right) ,$ $\mathrm{(i)}\Longleftrightarrow \mathrm{(ii)}%
\Longleftrightarrow \mathrm{(iii)}\Longleftrightarrow \mathrm{(iv)}$ follows
from Corollary \ref{Cor.4.2}. The equivalence between $\alpha =\max (\mathrm{%
D}_{\mathcal{H}})$ and $\alpha =\max (\mathrm{D})$ under the condition on $%
\mathcal{H}$ comes from Proposition \ref{pro1}.\hfill $\square $

\begin{remark}
\label{ddd}Given a bounded LSIP problem $(\mathrm{P})$ with $x^{\ast }\neq
0_{n}$ and optimal value $\alpha ,$ it has been proved in \cite[Theorem 8.3]%
{GL98} that $\left( \mathrm{P}\right) $ is reducible (i.e., $\mathcal{F}(T)$%
-reducible) if and only if
\begin{equation}
-\left( x^{\ast },\alpha \right) \in \mathcal{K}\cup \left( \left\{
0_{n}\right\} \times \mathbb{R}_{+}\right) ,  \label{4.7}
\end{equation}%
with $\mathcal{K}\cup \left( \left\{ 0_{n}\right\} \times \mathbb{R}%
_{+}\right) $ being the so-called characteristic cone of the system $%
\{\left\langle a_{t}^{\ast },x\right\rangle \leq b_{t},\ \forall t\in T\}$
(e.g., \cite[Chap. 4]{GL98}). We now show that we can eliminate the
half-line $\left\{ 0_{n}\right\} \times \mathbb{R}_{+}$\ in (\ref{4.7}). Let
$\lambda \in \mathbb{R}_{+}^{\left( T\right) }$ and $\mu \geq 0$ be such
that
\begin{equation*}
-\left( x^{\ast },\alpha \right) =\sum\nolimits_{t\in T}\lambda
_{t}(a_{t}^{\ast },b_{t})+\mu \left( 0_{n},1\right) .
\end{equation*}%
Since $x^{\ast }\neq 0_{n},$ $\gamma :=\sum\nolimits_{t\in T}\lambda _{t}>0$
and we can write%
\begin{equation*}
-\left( x^{\ast },\alpha \right) =\sum\nolimits_{t\in T}\lambda _{t}\left[
(a_{t}^{\ast },b_{t})+\left( 0_{n},\frac{\mu }{\gamma }\right) \right] \in
\mathcal{K}.
\end{equation*}%
\newline
The reducibility of $(\mathrm{D})$ is not discussed in \cite{GL98}.
\end{remark}

\begin{example}
\label{Exam4.3} We now revisit Examples \ref{Exam4.1} and \ref{Exam4.2} in
the light of Corollary \ref{Cor.4.3}. For problem $\left( \mathrm{P}\right) $
in Example \ref{Exam4.1} one has
\begin{eqnarray*}
\mathcal{K}_{\mathcal{H}_{1}} &=&\tbigcup\limits_{t\in T}\mathbb{R}%
_{+}\left( \left\{ (a_{t}^{\ast },b_{t})\right\} +\left\{ 0_{2}\right\}
\times \mathbb{R}_{+}\right) \\
&=&\mathbb{R}_{+}\left( \left\{ \left( a_{t}^{\ast },b_{t}\right) ,\ t\in
T\right\} +(0_{2},\mathbb{R}_{+})\right) \\
&=&\mathbb{R}_{+}\left( \left\{ \left( -t,t-1,t^{2}-t)\right) ,\ t\in \left[
0,1\right] \right\} +(0_{2},\mathbb{R}_{+})\right) ,
\end{eqnarray*}%
where $\left\{ \left( -t,t-1,t^{2}-t)\right) ,t\in \left[ 0,1\right]
\right\} $ is an arch of parabola contained in the plane $x_{1}+x_{2}=-1$
which connects the points $\left( 0,-1,0\right) $ and $\left( -1,0,0\right)
. $ The equation of the cone generated by the parabola $\left\{ \left(
-t,t-1,t^{2}-t)\right) ,t\in \mathbb{R}\right\} $ is $x_{3}=\frac{x_{1}x_{2}%
}{x_{1}+x_{2}}.$ Defining $\psi :\mathbb{R}^{2}\longrightarrow \overline{%
\mathbb{R}}$ such that
\begin{equation*}
\psi \left( x\right) =\left\{
\begin{array}{ll}
\frac{x_{1}x_{2}}{x_{1}+x_{2}}, & x\in \mathbb{R}_{-}^{2}\diagdown \left\{
0_{2}\right\} , \\
0, & x=0_{2}, \\
+\infty , & \text{else,}%
\end{array}%
\right.
\end{equation*}%
as
\begin{equation*}
\nabla ^{2}\psi \left( x\right) =\frac{2}{\left( x_{1}+x_{2}\right) ^{3}}%
\left[
\begin{array}{cc}
-x_{2}^{2} & x_{1}x_{2} \\
x_{1}x_{2} & -x_{1}^{2}%
\end{array}%
\right] ,
\end{equation*}%
is positive semidefinite on the interior of $\func{dom}\psi =\mathbb{R}%
_{-}^{2},$ $\psi $ is convex and also on $\func{dom}\psi $ by continuity.
Since $\psi $ is also lsc, the convex cone $\mathcal{K}_{\mathcal{H}_{1}}=%
\limfunc{epi}\psi $ is closed. Since $\mathcal{H}_{1}$-strong duality holds
in the four cases discussed in Example \ref{Exam4.1},\ we conclude again,
applying now Corollary \ref{Cor.4.3}\textrm{(ii)}, with $\alpha =\frac{%
x_{1}^{\ast }x_{2}^{\ast }}{x_{1}^{\ast }+x_{2}^{\ast }}$ if $x^{\ast }\in
\mathbb{R}_{++}^{2}$ and $\alpha =0$ otherwise, that \textit{\ }$(\mathrm{P}%
) $ and $(\mathrm{D})$ are $\mathcal{H}_{1}$-reducible for all $x^{\ast }\in
\mathbb{R}_{+}^{2}.$\newline
Observe that the algebraic condition $-\left( x^{\ast },\alpha \right) \in
\mathcal{K}_{\mathcal{H}_{1}}$ for $\mathcal{H}_{1}$-reducibility is easier
to check than the $\mathcal{H}_{1}$-strong duality. In fact, if $x^{\ast
}\in \mathbb{R}_{++}^{2},$ the unique solution of
\begin{equation*}
-\left( x_{1}^{\ast },x_{2}^{\ast },\frac{x_{1}^{\ast }x_{2}^{\ast }}{%
x_{1}^{\ast }+x_{2}^{\ast }}\right) =\lambda _{t}\left(
-t,t-1,t^{2}-t\right) ,\text{ }\lambda _{t}\geq 0,t\in \left[ 0,1\right] ,
\end{equation*}%
is $\left( t,\lambda _{t}\right) =\left( \frac{x_{1}^{\ast }}{x_{1}^{\ast
}+x_{2}^{\ast }},x_{1}^{\ast }+x_{2}^{\ast }\right) .$ Similarly, one has
that $-\left( x^{\ast },0\right) \in \mathcal{K}_{\mathcal{H}}$ holds when $%
x^{\ast }$ lies on the boundary of $\mathbb{R}_{+}^{2}.$ \newline
Regarding Example \ref{Exam4.2}, $\mathcal{K}_{\mathcal{H}_{1}}=\left\{
0_{n+1}\right\} \cup \limfunc{epi}\nolimits_{s}\psi $ is a convex but
non-closed cone and, since $(\mathrm{P})$ is not $\mathcal{H}_{1}$%
-reducible, the same happens with $(\mathrm{D})$\textit{. }Here the
condition $-\left( x^{\ast },\alpha \right) \in \mathcal{K}_{\mathcal{H}%
_{1}} $ writes%
\begin{equation*}
-\left( x^{\ast },\frac{x_{1}^{\ast }x_{2}^{\ast }}{x_{1}^{\ast
}+x_{2}^{\ast }}\right) =\lambda _{t}\left( -t,t-1,t^{2}-t+\frac{1}{r}%
\right) ,\text{ }\lambda _{t}\geq 0,t\in \left[ 0,1\right] ,r\in \mathbb{N},
\end{equation*}%
which has no solution for any $x^{\ast }\in \mathbb{R}_{+}^{2}.$
\end{example}

\section{$\mathcal{H}$-stable strong duality}

Let us go back to the general CIP problem \textrm{(P)} in (\ref{1.1}). By (%
\ref{2.1bis}) we have\textbf{\ }%
\begin{equation}
(f+\delta _{E})^{\ast }(x^{\ast })\leq \left( f+\sum_{t\in H}\mu
_{t}f_{t}\right) ^{\ast }(x^{\ast }),\ \text{for all }x^{\ast }\in X^{\ast
},\ H\in \mathcal{F}(T)\text{, and\ }\mu \in \mathbb{R}_{+}^{H}.  \label{51a}
\end{equation}

\begin{definition}
\label{def61} We say that the $\mathcal{H}$-strong\textit{\ }duality for (P)
holds at a given $x^{\ast }\in X^{\ast }$ if
\begin{equation}
(f+\delta _{E})^{\ast }(x^{\ast })=\min_{H\in \mathcal{H}\text{,\ }\mu \in
\mathbb{R}_{+}^{H}}\left( f+\sum_{t\in H}\mu _{t}f_{t}\right) ^{\ast
}(x^{\ast }).  \label{eq1}
\end{equation}
\end{definition}

For $x^{\ast }=0_{X^{\ast }}$, \eqref{eq1} amounts to the relation
\begin{equation}
\inf (\mathrm{P})=\max (\mathrm{D}_{\mathcal{H}}).  \label{marco51}
\end{equation}

\begin{definition}
If \eqref{eq1} is satisfied for all $x^{\ast }\in X^{\ast }$, one says that
\textit{$\mathcal{H}$-stable strong duality} for \textrm{(P)} {holds}, which
amounts to say that
\begin{equation*}
\limfunc{epi}(f+\delta _{E})^{\ast }=\mathcal{A}_{\mathcal{H}}.
\end{equation*}
\end{definition}

We first consider the $\mathcal{H}$-strong duality in the general case. We
first recall some general facts whose elementary proofs are similar to those
in \cite{Volle2015} and use (\ref{3.3k}).

\begin{lemma}
\label{lem51aa} The following assertions always hold:\newline
$\mathrm{(i)}$ $\mathrm{epi}_{s}\varphi _{\mathcal{H}}\subset \mathcal{A}_{%
\mathcal{H}}\subset \mathrm{epi}\varphi _{\mathcal{H}}\subset \mathrm{epi}%
(f+\delta _{E})^{\ast }$.\newline
$\mathrm{(ii)}$ $\overline{\mathrm{co}}\left( \mathrm{epi}\varphi _{\mathcal{%
H}}\right) =\overline{\mathrm{co}}\left( \mathcal{A}_{\mathcal{H}}\right)
\subset \mathrm{epi}\left( f+\delta _{E}\right) ^{\ast }$.
\end{lemma}

Another lemma will be useful.

\begin{lemma}
\label{lem1} Let $\mathcal{H}\subset \mathcal{F}(T)$ be a covering family.
Assume that the convex proper functions $\{f;\ f_{t},t\in T\}$ are lsc (in
other words $\{f;\ f_{t},t\in T\}\subset \Gamma (X)$). Then, we have
\begin{equation*}
(\varphi _{\mathcal{H}})^{\ast }=f+\delta _{E}.
\end{equation*}%
If, moreover, $E\cap (\func{dom}f)\neq \emptyset $, then $\limfunc{epi}%
(f+\delta _{E})^{\ast }$ coincides with the $w^{\ast }$-closed convex hull
of $\mathcal{A}_{\mathcal{H}}$, namely,
\begin{equation*}
\limfunc{epi}(f+\delta _{E})^{\ast }=\overline{\limfunc{co}}\left(
\bigcup\limits_{_{\substack{ H\in \mathcal{H}  \\ \mu \in \mathbb{R}_{+}^{H}
}}}\limfunc{epi}\left( f+\sum_{t\in H}\mu _{t}f_{t}\right) ^{\ast }\right) .
\end{equation*}
\end{lemma}

\noindent \textbf{Proof }We have
\begin{align*}
(\varphi _{\mathcal{H}})^{\ast }& =\left( \inf\limits_{_{\substack{ H\in
\mathcal{H}  \\ \mu \in \mathbb{R}_{+}^{H}}}}\left( f+\sum_{t\in H}\mu
_{t}f_{t}\right) ^{\ast }\right) ^{\ast } \\
& =\sup\limits_{_{\substack{ H\in \mathcal{H}  \\ \mu \in \mathbb{R}_{+}^{H}
}}}\left( f+\sum_{t\in H}\mu _{t}f_{t}\right) ^{\ast \ast }=\sup\limits_{
_{\substack{ H\in \mathcal{H}  \\ \mu \in \mathbb{R}_{+}^{H}}}}\left(
f+\sum_{t\in H}\mu _{t}f_{t}\right) \\
& =f+\sup\limits_{_{\substack{ H\in \mathcal{H}  \\ \mu \in \mathbb{R}%
_{+}^{H} }}}\left( \sum_{t\in H}\mu _{t}f_{t}\right) =f+\sup\limits_{H\in
\mathcal{H}}\sup\limits_{\mu \in \mathbb{R}_{+}^{H}}\left( \sum_{t\in H}\mu
_{t}f_{t}\right) \\
& =f+\sup\limits_{H\in \mathcal{H}}\delta _{E_{H}}=f+\delta _{E_{\mathcal{H}%
}}=f+\delta _{E}.
\end{align*}%
If, moreover, $E\cap (\func{dom}f)\neq \emptyset $, then $\func{dom}(\varphi
_{\mathcal{H}})^{\ast }\neq \emptyset $ and by \cite[Proposition 3.2]{ET74}%
,\smallskip\
\begin{equation*}
\limfunc{epi}(f+\delta _{E})^{\ast }=\limfunc{epi}(\varphi _{\mathcal{H}%
})^{\ast \ast }=\overline{\limfunc{co}}(\limfunc{epi}\varphi _{\mathcal{H}%
}).\smallskip
\end{equation*}%
Lemma \ref{lem51aa}(ii) concludes the proof.\hfill $\square $

\begin{theorem}[$\mathcal{H}$-stable strong duality]
\label{thm2} Let $\mathcal{H}\subset \mathcal{F}(T)$ be a covering family, $%
x^{\ast }\in X^{\ast }$, and consider the following statements:\newline
$\left( \mathrm{i}\right) $ $\mathcal{H}$-strong duality holds at $x^{\ast }$%
.\newline
$\left( \mathrm{ii}\right) $ $\mathcal{A}_{\mathcal{H}}$ is $w^{\ast }$%
-closed convex regarding $\{x^{\ast }\}\times \mathbb{R}$.

Then we have $\mathrm{(i)\Rightarrow (ii)}$. If, moreover, $\{f;\ f_{t},t\in
T\}\subset \Gamma (X)$ and $E\cap \func{dom}f\neq \emptyset $, then $\mathrm{%
(i)\Leftrightarrow (ii)}$. In particular, $\mathcal{H}$-stable strong
duality holds if and only if $\mathcal{A}_{\mathcal{H}}$ is $w^{\ast }$%
-closed convex.
\end{theorem}

\noindent \textbf{Proof }$\left[ \mathrm{(i)\Rightarrow (ii)}\right] $ Let $%
(x^{\ast },r)\in \overline{\limfunc{co}}(\mathcal{A}_{\mathcal{H}})$. By
Lemma \ref{lem51aa}$\mathrm{(ii)}$ we have $(f+\delta _{E})^{\ast }(x^{\ast
})\leq r$ and, by assumption \textrm{(i),} there exist $H\in \mathcal{H}$, $%
\mu \in \mathbb{R}_{+}^{H}$ such that $\left( f+\delta _{E}\right) ^{\ast
}(x^{\ast })=\left( f+\sum_{t\in H}\mu _{t}f_{t}\right) ^{\ast }(x^{\ast
})\leq r$, that means $(x^{\ast },r)\in \mathcal{A}_{H}\subset \mathcal{A}_{%
\mathcal{H}}$. Therefore, $\mathcal{A}_{\mathcal{H}}$ is $w^{\ast }$-closed
regarding $\{x^{\ast }\}\times \mathbb{R}$.

$\left[ \mathrm{(ii)\Rightarrow (i)}\right] $ Since\textbf{\ }$E\cap \func{%
dom}f\neq \emptyset $, we have $(f+\delta _{E})^{\ast }(x^{\ast })\neq
-\infty $. If $(f+\delta _{E})^{\ast }(x^{\ast })=+\infty $, then for each $%
H\in \mathcal{H}$ and each $\mu \in \mathbb{R}_{+}^{H}$, we get, from %
\eqref{51a}, $\left( f+\sum_{t\in H}\mu _{t}f_{t}\right) ^{\ast }(x^{\ast
})=+\infty $ and
\begin{equation*}
(f+\delta _{E})^{\ast }(x^{\ast })=+\infty =\min_{H\in \mathcal{H},\,\mu \in
\mathbb{R}_{+}^{H}}\left( f+\sum_{t\in H}\mu _{t}f_{t}\right) ^{\ast
}(x^{\ast }).
\end{equation*}

Assume now that $r:=(f+\delta _{E})^{\ast }(x^{\ast })\in \mathbb{R}$. We
have $(x^{\ast },r)\in \limfunc{epi}(f+\delta _{E})^{\ast }=\overline{%
\limfunc{co}}(\mathcal{A}_{\mathcal{H}})$ (see Lemma \ref{lem1}). By
assumption \textrm{(ii)} we obtain that $(x^{\ast },r)\in \mathcal{A}_{%
\mathcal{H}}$, and there exist $\tilde{H}\in \mathcal{H}$, $\tilde{\mu}\in
\mathbb{R}_{+}^{\tilde{H}}$, such that
\begin{align*}
\left( f+\sum_{t\in \tilde{H}}\tilde{\mu}_{t}f_{t}\right) ^{\ast }(x^{\ast
})& \leq r=(f+\delta _{E})^{\ast }(x^{\ast }) \\
& \leq \inf\limits_{_{\substack{ H\in \mathcal{H}  \\ \mu \in \mathbb{R}%
_{+}^{H}}}}\left( f+\sum_{t\in H}\mu _{t}f_{t}\right) ^{\ast }(x^{\ast }) \\
& \leq \left( f+\sum_{t\in \tilde{H}}\tilde{\mu}_{t}f_{t}\right) ^{\ast
}(x^{\ast }),
\end{align*}%
and we are done.\hfill $\square \medskip $

Once again, regarding the problems in Examples \ref{Exam4.1} and \ref%
{Exam4.2}, we conclude from Theorem \ref{thm2} that $(\mathrm{P})\ $%
satisfies $\mathcal{H}_{1}$\textit{-}stable strong duality in the first case
because $\mathcal{A}_{\mathcal{H}}$ is $w^{\ast }$-closed convex, but not in
the second because $\mathcal{A}_{\mathcal{H}}$ is just convex.

\begin{lemma}
\label{lem2} For each family $\mathcal{H}\subset \mathcal{F}(T)$ we have
\begin{equation*}
\mathcal{A}_{\mathcal{H}}\text{ convex }\Rightarrow \varphi _{\mathcal{H}}%
\text{ convex }\Rightarrow \overline{\varphi _{\mathcal{H}}}\text{ convex }%
\Leftrightarrow \overline{\mathcal{A}_{\mathcal{H}}}\text{ convex }%
\Leftrightarrow \overline{\limfunc{co}}(\mathcal{A}_{\mathcal{H}})=\overline{%
\mathcal{A}_{\mathcal{H}}}.
\end{equation*}
\end{lemma}

\noindent \textbf{Proof }Since $\limfunc{epi}_{s}\varphi _{\mathcal{H}%
}\subset \mathcal{A}_{\mathcal{H}}\subset \limfunc{epi}\varphi _{\mathcal{H}%
} $ we have $\varphi _{\mathcal{H}}(x^{\ast })=\inf \{r\in \mathbb{R}%
:(x^{\ast },r)\in \mathcal{A}_{\mathcal{H}}\}$ and then,
\begin{eqnarray*}
\mathcal{A}_{\mathcal{H}}\text{ convex }\Rightarrow \varphi _{\mathcal{H}}%
\text{ convex }\Rightarrow \overline{\varphi _{\mathcal{H}}}\text{ convex }
&\Leftrightarrow &\overline{\limfunc{epi}\varphi _{\mathcal{H}}}=\overline{%
\mathcal{A}_{\mathcal{H}}}\text{ convex } \\
&\Rightarrow &\overline{\limfunc{co}}(\mathcal{A}_{\mathcal{H}})=\overline{%
\limfunc{co}}(\overline{\mathcal{A}_{\mathcal{H}}})=\overline{\mathcal{A}_{%
\mathcal{H}}}\Rightarrow \overline{\mathcal{A}_{\mathcal{H}}}\text{ convex}.
\end{eqnarray*}%
$\square $ %
%
%
%
%
%
%
%
%
%
%
%
%
%
%
%
%
%

As immediate consequence of Theorem \ref{thm2} and Lemma \ref{lem2} we
establish the following consequences:

\begin{corollary}
\label{cor1} Let $\mathcal{H}\subset \mathcal{F}(T)$ be a covering family
with $\overline{\mathcal{A}_{\mathcal{H}}}$ convex. Assume that $\{f;\
f_{t},t\in T\}\subset \Gamma (X)$ and $E\cap (\func{dom}f)\neq \emptyset $.
Then $\mathcal{H}$-strong duality holds at a given $x^{\ast }\in X^{\ast }$
if and only if $\mathcal{A}_{\mathcal{H}}$ is $w^{\ast }$-closed regarding $%
\{x^{\ast }\}\times \mathbb{R}$.

In particular, $\mathcal{H}$-stable strong duality holds if and only if $%
\mathcal{A}_{\mathcal{H}}$ is $w^{\ast }$-closed.
\end{corollary}

\begin{corollary}
\label{cor3} Assume that $\{f;\ f_{t},t\in T\}\subset \Gamma (X)$, $E\cap (%
\func{dom}f)\neq \emptyset $, and let $\mathcal{H}\subset \mathcal{F}(T)$ be
a directed covering family. Then $\mathcal{H}$-strong duality holds at a
given $x^{\ast }\in X^{\ast }$ if and only if $\mathcal{A}_{\mathcal{H}}$ is
$w^{\ast }$-closed regarding $\{x^{\ast }\}\times \mathbb{R}$. In
particular, $\mathcal{H}$-stable strong duality holds if and only if $%
\mathcal{A}_{\mathcal{H}}$ (alias $\mathcal{A}$) is $w^{\ast }$-closed.
\end{corollary}

\noindent \textbf{Proof }Since $\mathcal{A}_{\mathcal{H}}=\mathcal{A}$ (see
Proposition \ref{pro1}) and $\mathcal{A}$ is convex (recall (\ref{marco7})),
Corollary \ref{cor1} concludes the proof. \hfill $\square $

We now give a corollary, addressing the LIP problem in (\ref{LIP}), whose
proof is a straightforward consequence of Theorem \ref{thm2} and the
relation $\mathcal{A}_{\mathcal{H}}=\left( x^{\ast },0\right) +\mathcal{K}_{%
\mathcal{H}}$.

\begin{corollary}
\label{cor3bis} Consider the LIP problem $(\mathrm{P})$ in (\ref{LIP}), and
let $\mathcal{H}$ be a covering family. Then, the following statements are
equivalent:\newline
$(\mathrm{i})$ $\inf (\mathrm{P})=\max (\mathrm{D}_{\mathcal{H}}).$\newline
$(\mathrm{ii})$ $\overline{\limfunc{co}}(\mathcal{K}_{\mathcal{H}})\cap
(\{-x^{\ast }\}\times \mathbb{R})=\mathcal{K}_{\mathcal{H}}\cap (\{-x^{\ast
}\}\times \mathbb{R}).$\newline
If $\mathcal{H}$ is additionally directed, $(\mathrm{i})$ is equivalent to%
\newline
$(\mathrm{iii})$ $\overline{\mathcal{K}_{\mathcal{H}}}\cap (\{-x^{\ast
}\}\times \mathbb{R})=\mathcal{K}_{\mathcal{H}}\cap (\{-x^{\ast }\}\times
\mathbb{R}),$\newline
and $\mathcal{H}$-stable strong stability duality holds if and only if $%
\mathcal{K}_{\mathcal{H}}$ is w$^{\ast }$-closed.
\end{corollary}

Consider again problem $(\mathrm{P})$ in Example \ref{Exam4.1}, where we
proved that $\mathcal{K}_{\mathcal{H}_{1}}$ is closed and convex. So, we
conclude again, now from Corollary \ref{cor3bis}, that $\mathcal{H}_{1}$%
-strong duality holds.\medskip

We now give a new glimpse on $\mathcal{H}_{\mathbb{N}}$-stable strong
duality for convex infinite countable programs. Given the family $\{f;\
f_{k},k\in \mathbb{N}\}\subset \Gamma (X)$, consider the countable convex
optimization problem
\begin{equation*}
(\mathrm{P}_{\mathbb{N}})\quad \inf f(x)\,\,\mathrm{s.t.}\,f_{k}(x)\leq
0,k\in \mathbb{N},
\end{equation*}%
and, for each $m\in \mathbb{N}$, the finite subproblem
\begin{equation*}
(\mathrm{P}_{m})\quad \inf f(x)\,\,\mathrm{s.t.}\,f_{k}(x)\leq 0,\ k\in
\{1,\cdots ,m\}.
\end{equation*}%
The Lagrangian-Haar dual of $(\mathrm{P}_{\mathbb{N}})$ reads
\begin{equation*}
(\mathrm{D}_{\mathbb{N}})\quad \sup_{\lambda \in \mathbb{R}_{+}^{(\mathbb{N}%
)}}\inf_{x\in X}\left\{ f(x)+\sum_{k\in \mathbb{N}}\lambda
_{k}f_{k}(x)\right\} .
\end{equation*}%
Let us consider the Lagrangian dual of the subproblem $(\mathrm{P}_{m})$,
that is
\begin{equation*}
(\mathrm{D}_{m})\quad \sup_{\mu \in \mathbb{R}_{+}^{m}}\inf_{x\in X}\left\{
f(x)+\sum_{k=1}^{m}\mu _{k}f_{k}(x)\right\} .
\end{equation*}%
In terms of $\mathcal{H}$-duality, the corresponding family is
\begin{equation*}
\mathcal{H}_{\mathbb{N}}=\left\{ \{1,\cdots ,m\}\right\} _{m\in \mathbb{N}},
\end{equation*}%
which is a directed covering family. Then, by Proposition \ref{pro1}, we have%
\begin{equation}
\mathcal{A}=\bigcup_{\lambda \in \mathbb{R}_{+}^{(\mathbb{N})}}\limfunc{epi}%
\left( f+\sum_{k\in \mathbb{N}}\lambda _{k}f_{k}\right) ^{\ast
}=\bigcup_{m\in \mathbb{N},\mu \in \mathbb{R}_{+}^{m}}\limfunc{epi}\left(
f+\sum_{k=1}^{m}\mu _{k}f_{k}\right) ^{\ast }=\mathcal{A}_{\mathcal{H}_{%
\mathbb{N}}},  \label{eqn2}
\end{equation}%
which is a convex subset of $X^{\ast }\times \mathbb{R}$, and also
\begin{equation}
\sup (\mathrm{D}_{\mathbb{N}})=\sup_{m\in \mathbb{N}}\sup (\mathrm{D}%
_{m})=\lim_{m\rightarrow \infty }\sup (\mathrm{D}_{m}),  \label{marco43}
\end{equation}%
showing that the optimal value of $(\mathrm{D}_{\mathbb{N}})$ can be
arbitrarily approached by solving a sequence of finite subproblems.

\begin{theorem}[$\mathcal{H}_{\mathbb{N}}$-stable strong duality]
\label{thm3} Assume that $E\cap (\func{dom}f)\neq \varnothing $ and let $%
x^{\ast }\in X^{\ast }$. The following statements are equivalent:\newline
$\left( \mathrm{i}\right) $ $(f+\delta _{E})^{\ast }(x^{\ast
})=\min\limits_{\mu \in \mathbb{R}_{+}^{m}}\left( f+\sum\limits_{k=1}^{m}\mu
_{k}f_{k}\right) ^{\ast }(x^{\ast })$ for all $m\in \mathbb{N}$ sufficiently
large.\newline
$\left( \mathrm{ii}\right) $ $\mathcal{A}$ is $w^{\ast }$-closed regarding $%
\{x^{\ast }\}\times \mathbb{R}$.
\end{theorem}

\noindent \textbf{Proof }$\left[ \mathrm{(i)\ \Rightarrow (ii)}\right] $
There exist $m\in \mathbb{N}$ and $\mu \in \mathbb{R}_{+}^{m}$ such that $%
(f+\delta _{E})^{\ast }(x^{\ast })=\left( f+\sum_{k=1}^{m}\mu
_{k}f_{k}\right) ^{\ast }(x^{\ast })$. Consequently, $\mathcal{H}_{\mathbb{N}%
}$-strong duality holds at $x^{\ast }$. By Corollary \ref{cor3}, $\mathcal{A}
$ is $w^{\ast }$-closed regarding $\{x^{\ast }\}\times \mathbb{R}$.\newline

$\left[ \mathrm{(ii)\ \Rightarrow (i)}\right] $ By Corollary \ref{cor3} and
the equality $\mathcal{A}=\mathcal{A}_{\mathcal{H}_{\mathbb{N}}}$, there
exist $N\in \mathbb{N}$ and $\bar{\mu}\in \mathbb{R}_{+}^{N}$ such that $%
(f+\delta _{E})^{\ast }(x^{\ast })=\left( f+\sum_{k=1}^{N}\bar{\mu}%
_{k}f_{k}\right) ^{\ast }(x^{\ast })$. For each $m\geq N$, let us now define
\begin{equation*}
\tilde{\mu}_{k}=%
\begin{cases}
\bar{\mu}_{k}, & \text{ if }1\leq k\leq N, \\
0, & \text{ if }N<k\leq m.%
\end{cases}%
\end{equation*}%
We have $\tilde{\mu}\in \mathbb{R}_{+}^{m}$ and $\sum_{k=1}^{m}\widetilde{%
\mu }_{k}f_{k}=\sum_{k=1}^{N}\bar{\mu}_{k}f_{k}$. Finally, if $m\geq N,$ we
have
\begin{eqnarray*}
\inf_{m\geq N}\inf_{\mu \in \mathbb{R}_{+}^{p}}\left( f+\sum_{k=1}^{m}\mu
_{k}f_{k}\right) ^{\ast }(x^{\ast }) &\leq &\left( f+\sum_{k=1}^{m}\tilde{\mu%
}_{k}f_{k}\right) ^{\ast }(x^{\ast }) \\
&=&(f+\delta _{E})^{\ast }(x^{\ast })\leq \inf_{p\in \mathbb{N},\ \mu \in
\mathbb{R}_{+}^{p}}\left( f+\sum_{k=1}^{p}\mu _{k}f_{k}\right) ^{\ast
}(x^{\ast }),
\end{eqnarray*}%
and \textrm{(i)} holds.\hfill $\square $

\begin{corollary}
\label{cor4} Assume that $\inf (\mathrm{P}_{\mathbb{N}})\neq +\infty $. The
following statements are equivalent:\newline
$\left( \mathrm{i}\right) $ $\inf (\mathrm{P}_{\mathbb{N}})=\max (\mathrm{D}%
_{m})$ for all $m$ sufficiently large.\newline
$\left( \mathrm{ii}\right) $ $\mathcal{A}$ is $w^{\ast }$-closed regarding $%
\{0_{X^{\ast }}\}\times \mathbb{R}$.
\end{corollary}

\noindent \textbf{Proof }Apply Theorem \ref{thm3} with $x^{\ast }=0_{X^{\ast
}}$.\hfill $\square \medskip $

\begin{remark}
\label{rem55} For each $m\in \mathbb{N}$ it holds that
\begin{equation*}
\inf (\mathrm{P}_{\mathbb{N}})\ \geq \
\begin{array}{c}
\inf (\mathrm{P}_{m}) \\
\sup (\mathrm{D}_{\mathbb{N}})%
\end{array}%
\ \geq \ \ \sup (\mathrm{D}_{m}).
\end{equation*}%
Consequently, condition $(\mathrm{ii})$ of Corollary \ref{cor4} guarantees
that, for $m$ sufficient large,
\begin{equation*}
\inf (\mathrm{P}_{\mathbb{N}})\ =\ \inf (\mathrm{P}_{m})\ =\ \sup (\mathrm{D}%
_{m})\ =\ \ \sup (\mathrm{D}_{\mathbb{N}}),
\end{equation*}%
and, in particular, that $(\mathrm{P}_{\mathbb{N}})$ and $(\mathrm{D}_{%
\mathbb{N}})$ are simultaneously $\mathcal{H}_{\mathbb{N}}$-reducible. Note
that for convex SIP problems such that $\{f;\ f_{t},t\in T\}\subset \Gamma (%
\mathbb{R}^{n})$ one can find in \cite[Theorem 4.1 and Corollary 4.2]%
{Karney1983} a recession condition ensuring that $\inf (\mathrm{P}_{m})=\sup
(\mathrm{D}_{m})\ $for $m$ sufficiently large, with $\inf (\mathrm{P}%
_{m})\rightarrow \inf (\mathrm{P}_{\mathbb{N}})$ as $m\rightarrow \infty .$
\end{remark}

We finish this section providing easily checkable conditions guaranteeing the%
\textbf{\ }$\mathcal{H}_{1}$-(stable) strong duality for \textrm{(P)}.
Recall that the $\mathcal{H}_{1}$-dual of $(\mathrm{P})$ reads
\begin{equation*}
(\mathrm{D}_{\mathcal{H}_{1}})\ \ \ \ \ \ \sup\limits_{(t,\mu )\in T\times
\mathbb{R}_{+}}\inf_{x\in X}\left\{ f(x)+\mu f_{t}(x)\right\} .
\end{equation*}

\begin{theorem}[$\mathcal{H}_{1}$-stable strong duality]
\label{thm1} Assume:\newline
$\left( \mathrm{a}\right) $ $\func{dom}f\subset \bigcap\limits_{t\in T}\func{%
dom}f_{t}$.\newline
$\left( \mathrm{b}\right) $ $T$ is a convex and compact subset of some
locally convex space.\newline
$\left( \mathrm{c}\right) $ $T\ni t\mapsto f_{t}(x)$ is concave and usc on $%
T $ for each $x\in \bigcap\nolimits_{t\in T}\func{dom}f_{t}.$\newline
$\left( \mathrm{d}\right) $\textrm{\ }There exists $\bar{x}\in \func{dom}f$
such that $f_{t}(\bar{x})<0$ for all $t\in T$.\newline
Then $\mathcal{H}_{1}$-stable strong duality holds. In particular, we have
\begin{equation}
-\infty \leq \inf (\mathrm{P})=\max (\mathrm{D}_{\mathcal{H}_{1}})<+\infty .
\label{4.1}
\end{equation}
\end{theorem}

\noindent \textbf{Proof }We first prove (\ref{4.1}). Let us consider the
convex function $h:=\sup_{t\in T}f_{t}$. Thanks to $\left( \mathrm{d}\right)
$ we have\
\begin{equation*}
\infty \leq \inf (\mathrm{P})=\inf \{f(x):x\in \lbrack h\leq 0]\cap \func{dom%
}f\}<+\infty .
\end{equation*}%
The compactness of $T,$ assumed in $\left( \mathrm{b}\right) ,$ and the
upper semicontinuity of the functions $t\mapsto f_{t}(x)$ on $T,$ for each $%
x\in \bigcap\nolimits_{t\in T}\func{dom}f_{t},$ assumed in $\left( \mathrm{c}%
\right) $, yield
\begin{equation*}
\func{dom}h=\bigcap\limits_{t\in T}\func{dom}f_{t}\supseteq \func{dom}f\text{
and }\bar{x}\in \lbrack h<0]\cap \func{dom}f.
\end{equation*}%
Since\textbf{\ }$h$ is convex and proper, by \cite[Theorem 2.9.3]{Z02} there
exists $\bar{\mu}\in \mathbb{R}_{+}$ such that
\begin{equation*}
\inf (\mathrm{P})=\inf_{h(x)\leq 0}f(x)=\inf_{x\in \func{dom}f}\{f(x)+\bar{%
\mu}h(x)\}=\inf_{x\in \func{dom}f}\max_{t\in T}\{f(x)+\bar{\mu}f_{t}(x)\}.
\end{equation*}%
Now, by the general minimax theorem \cite[Theorem 2.10.2]{Z02}, we have
\begin{equation*}
\inf (\mathrm{P})=\max_{t\in T}\inf_{x\in \func{dom}f}\{f(x)+\bar{\mu}%
f_{t}(x)\}.
\end{equation*}%
Finally, there exists $(\bar{t},\bar{\mu})\in T\times \mathbb{R}_{+}$ such
that
\begin{equation*}
\inf (\mathrm{P})=\inf_{x\in \func{dom}f}\{f(x)+\bar{\mu}f_{\bar{t}%
}(x)\}\leq \sup (\mathrm{D}_{\mathcal{H}_{1}})=\sup (\mathrm{D})\leq \inf (%
\mathrm{P}),
\end{equation*}%
which shows that (\ref{4.1}) holds.

Now, given an arbitrary $x^{\ast }\in X^{\ast },$ we can apply (\ref{4.1}),
replacing $f$ by $f-x^{\ast }$, since the corresponding assumptions \textrm{%
(a)}, \textrm{(b)}, \textrm{(c)}, \textrm{(d)} are the same. So, $\mathcal{H}%
_{1}$-stable strong duality for \textrm{(P)} holds.\hfill $\square \medskip $

Consider the problems in Examples \ref{Exam4.1} and \ref{Exam4.2}, with a
fixed objective functional $c^{\ast }$ instead of $x^{\ast }.$ The problem
in Example \ref{Exam4.1} enjoys $\mathcal{H}_{1}$\textit{-}stable strong
duality by Theorem \ref{thm1}, whose four assumptions trivially hold.
However, we cannot apply Theorem \ref{thm1} to the problem in Example \ref%
{Exam4.3} even though\ \textrm{(a)}, \textrm{(c)}, and \textrm{(d)} hold (in
fact, we have seen in different ways that $\mathcal{H}_{1}$\textit{-}stable
strong duality fails).

\begin{remark}
\label{rem2}\cite[Theorem 2.9.3]{Z02} is established under the rule $0\times
(+\infty )=+\infty $ (see \cite[p.39]{Z02}) instead of the rule $0\times
(+\infty )=0$ we use in this paper. In fact, since $\func{dom}f\subset \func{%
dom}h$, the relations $\inf_{h(x)\leq 0}f(x)=\inf_{x\in X}\{f(x)+\bar{\mu}%
h(x)\}=\inf_{x\in \func{dom}f}\{f(x)+\bar{\mu}h(x)\}$ we use in the proof of
Theorem \ref{thm1} are valid with both rules (recall that $%
h:=\sup\nolimits_{t\in T}f_{t}$).
\end{remark}

\begin{remark}
\label{rem3} In the proof of Theorem \ref{thm1} we use two fundamental
formulas. The first one, by \cite[Theorem 2.9.3]{Z02}, says that
\begin{equation}
\inf_{h(x)\leq 0}f(x)=\max_{\mu \geq 0}\inf_{x\in \func{dom}f}\{f(x)+\mu
h(x)\},\newline
\label{5.1}
\end{equation}%
in which the linear space $X\supset \func{dom}f$\ is in fact a locally
convex Hausdorff topological vector space. However, (\ref{5.1}) holds in
general linear spaces $X$ provided Slater condition
\begin{equation}
\exists \overline{x}\in \func{dom}f:h(\overline{x})<0  \label{5.2}
\end{equation}%
is satisfied. This can be found for instance in \cite[Lemma 1]{MV98}.\newline
The second formula says that, for $\overline{\mu }\geq 0,$%
\begin{equation*}
\inf_{x\in \func{dom}f}\max_{t\in T}\{f(x)+\bar{\mu}f_{t}(x)\}=\max_{t\in
T}\inf_{x\in \func{dom}f}\{f(x)+\bar{\mu}f_{t}(x)\},
\end{equation*}%
which is a direct consequence of \cite[Theorem 2.10.2]{Z02} with $X$ being
an arbitrary linear space and $T$ a compact convex set. By $\left( \mathrm{a}%
\right) $ and $\left( \mathrm{c}\right) ,$ $\left( \mathrm{d}\right) $ is
equivalent to (\ref{5.2}). Hence, Theorem \ref{thm1} remains valid if $X$ is
only required to be a linear space.\ \newline
Observe also that the last part of the argument remains true replacing $%
x^{\ast }\in X^{\ast }$\ by an arbitrary algebraic linear form $\ell $ on $X$
since assumptions $\left( \mathrm{a}\right) $, $\left( \mathrm{b}\right) $, $%
\left( \mathrm{c}\right) $, $\left( \mathrm{d}\right) $ are the same for $f$
and for $f-\ell $.
\end{remark}

\begin{remark}
Under the assumptions of Theorem \ref{thm1} we know that $\mathcal{H}_{1}$%
-stable strong duality holds (see Remark \ref{rem3}). By the first part of
Theorem \ref{thm2} we obtain that $\mathcal{A}_{\mathcal{H}%
_{1}}=\bigcup_{(t,\mu )\in T\times \mathbb{R}_{+}}\limfunc{epi}(f+\mu
f_{t})^{\ast }$ is $w^{\ast }$-closed convex.
\end{remark}

\begin{remark}
\label{rem4} If one retains only the assumptions $\left( \mathrm{a}\right) $%
, $\left( \mathrm{b}\right) $, $\left( \mathrm{c}\right) $ of Theorem \ref%
{thm1} (dropping the Slater condition $\left( \mathrm{d}\right) $) we have
that $\varphi _{\mathcal{H}_{1}}$ is convex. Let us prove this fact that
will be needed further on. Given $x^{\ast }\in X^{\ast },$ we have
\begin{equation*}
\varphi _{\mathcal{H}_{1}}(x^{\ast })=\inf_{\mu \in \mathbb{R}%
_{+}}\inf_{t\in T}\sup_{x\in \func{dom}f}\left\{ \langle x^{\ast },x\rangle
-f(x)-\mu f_{t}(x)\right\} .
\end{equation*}%
Applying again \cite[Theorem 2.10.2]{Z02} we have, for each $\mu \geq 0$,
\begin{align*}
\inf_{t\in T}\sup_{x\in \func{dom}f}\left\{ \langle x^{\ast },x\rangle
-f(x)-\mu f_{t}(x)\right\} & =\sup_{x\in \func{dom}f}\left\{ \langle x^{\ast
},x\rangle -f(x)+\min_{t\in T}(-\mu f_{t}(x))\right\} \\
& =\sup_{x\in \func{dom}f}\left\{ \langle x^{\ast },x\rangle -f(x)-\mu
\max_{t\in T}f_{t}(x)\right\} \\
& =(f+\mu h)^{\ast }(x^{\ast }).
\end{align*}%
Note that, for each $x\in \func{dom}f,$ the function
\begin{equation*}
\mathbb{R}_{+}\times X^{\ast }\ni (\mu ,x^{\ast })\mapsto \langle x^{\ast
},x\rangle -f(x)-\mu h(x)
\end{equation*}%
is affine. Consequently, the function
\begin{equation*}
\mathbb{R}_{+}\times X^{\ast }\ni (\mu ,x^{\ast })\mapsto (f+\mu h)^{\ast
}(x^{\ast })
\end{equation*}%
is convex and, finally, $\varphi _{\mathcal{H}_{1}}(x^{\ast })=\inf_{\mu
\geq 0}(f+\mu h)^{\ast }(x^{\ast })$ is convex too (e.g., \cite[Theorem
2.1.3(v)]{Z02}).\medskip
\end{remark}

\begin{corollary}
\label{cor2} Assume that $\{f;\ f_{t},t\in T\}\subset \Gamma (X)$, $E\cap
\func{dom}f\neq \emptyset $, and the conditions $\left( \mathrm{a}\right) $,
$\left( \mathrm{b}\right) $, $\left( \mathrm{c}\right) $ in Theorem \ref%
{thm1} are satisfied. Then, $\mathcal{H}_{1}$-stable strong duality holds if
and only if $\bigcup\limits_{(t,\mu )\in T\times \mathbb{R}_{+}}\limfunc{epi}%
(f+\mu f_{t})^{\ast }$ is $w^{\ast }$-closed.
\end{corollary}

\noindent \textbf{Proof }As mentioned in Remark \ref{rem4}, $\varphi _{%
\mathcal{H}_{1}}$ is convex. Consequently, by Lemma \ref{lem2}, $\overline{%
\mathcal{A}_{\mathcal{H}_{1}}}$ is convex and conclusion follows from
Corollary \ref{cor1}. \hfill $\square \medskip $

\section{Zero $\mathcal{H}$-duality gap}

In this section, we consider the general CIP problem in (\ref{1.1}) with the
feasible set $E=\cap _{t\in T}[f_{t}\leq 0]$. Given $H\in \mathcal{F}(T)$,
and $\emptyset \not=\mathcal{H}\subset \mathcal{F}(T)$, recall the sets $%
\mathcal{A}_{H}$ and $\mathcal{A}_{\mathcal{H}}$ defined in (\ref{2.2b}) and
(\ref{2.2bb}), respectively, and the function $\varphi _{\mathcal{H}}$ in (%
\ref{3.3b}) as follows:
\begin{eqnarray*}
\mathcal{A}_{H} &=&\bigcup\limits_{\mu \in \mathbb{R}_{+}^{H}}\limfunc{epi}%
\left( f+\sum_{t\in H}\mu _{t}f_{t}\right) ^{\ast },\ \ \ \ \mathcal{A}_{%
\mathcal{H}}=\bigcup\limits_{H\in \mathcal{H}}\mathcal{A}_{H}, \\
\varphi _{\mathcal{H}} &=&\inf_{H\in \mathcal{H}}\inf_{\mu \in \mathbb{R}%
_{+}^{H}}\left( f+\sum_{t\in H}\mu _{t}f_{t}\right) ^{\ast
}=\inf\limits_{H\in \mathcal{H}}\varphi _{H}.
\end{eqnarray*}

\begin{definition}
Given $\emptyset \not=\mathcal{H}\subset \mathcal{F}(T)$ and $x^{\ast }\in
X^{\ast }$, one says that $\mathcal{H}$-duality for $(\mathrm{P})$ holds at $%
x^{\ast }$ if
\begin{equation*}
\left( f+\delta _{E}\right) ^{\ast }(x^{\ast })=\varphi _{\mathcal{H}%
}(x^{\ast }).
\end{equation*}%
For $x^{\ast }=0_{X^{\ast }}$ that leads us to
\begin{equation*}
\inf \mathrm{(P)}=\sup (\mathrm{D}_{\mathcal{H}})=\sup\limits_{H\in \mathcal{%
H},\mu \in \mathbb{R}_{+}^{H}}\inf\limits_{x\in X}\left(
f(x)+\sum\limits_{t\in H}\mu _{t}f_{t}(x)\right) .
\end{equation*}
\end{definition}

We now characterize the $\mathcal{H}$-duality for $(\mathrm{P})$.

\begin{theorem}[Zero $\mathcal{H}$-duality gap]
\label{thm41a} Consider the following statements:\newline
$\left( \mathrm{i}\right) $ $\mathcal{H}$-duality for $(\mathrm{P})$ holds
at $x^{\ast }$, \newline
$\left( \mathrm{ii}\right) $ $\left( \overline{\mathrm{co}}\mathcal{A}_{%
\mathcal{H}}\right) \cap \left( \{x^{\ast }\}\times \mathbb{R}\right) =%
\overline{\mathcal{A}_{\mathcal{H}}\cap \left( \{x^{\ast }\}\times \mathbb{R}%
\right) }$.\newline
Then $\left( \mathrm{i}\right) \Rightarrow \left( \mathrm{ii}\right) $. If,
moreover, $\{f;$ $f_{t},t\in T\}\subset \Gamma (X)$, $E\cap \limfunc{dom}%
f\not=\emptyset $, and $\mathcal{H}$ is covering, then $\left( \mathrm{i}%
\right) \Leftrightarrow \left( \mathrm{ii}\right) $.
\end{theorem}

\noindent \textbf{Proof } $[\left( \mathrm{i}\right) \Rightarrow \left(
\mathrm{ii}\right) ]$ The inclusion $\left[ \supset \right] $ in $\left(
\mathrm{ii}\right) $ is obvious. Let us prove the opposite one. Let $%
(x^{\ast },r)\in \overline{\mathrm{co}}\mathcal{A}_{\mathcal{H}}$. We have
to check that $(x^{\ast },r)\in \overline{\mathcal{A}_{\mathcal{H}}\cap
\left( \{x^{\ast }\}\times \mathbb{R}\right) }$. By Lemma \ref{lem51aa}$%
\left( \mathrm{ii}\right) $, we have $(f+\delta _{E}(x^{\ast })\leq r$, and
by assumption $\left( \mathrm{i}\right) $, $\varphi _{\mathcal{H}}(x^{\ast
})\leq r<r+1/n$ for any $n\geq 1$.

On the other hand, it follows from Lemma \ref{lem51aa}(i) that
\begin{equation*}
(x^{\ast },r+1/n)\in \mathcal{A}_{\mathcal{H}}\cap \left( \{x^{\ast
}\}\times \mathbb{R}\right) ,\ \forall n\geq 1,
\end{equation*}%
and finally, $(x^{\ast },r)=\lim\limits_{n\rightarrow \infty }(x^{\ast
},r+1/n)\in \overline{\mathcal{A}_{\mathcal{H}}\cap \left( \{x^{\ast
}\}\times \mathbb{R}\right) }.$

We now assume that $\{f;$ $f_{t},t\in T\}\subset \Gamma (X)$, $E\cap
\limfunc{dom}\not=\emptyset $, $\mathcal{H}$ is covering, and we prove $%
\left( \mathrm{ii}\right) \Rightarrow \left( \mathrm{i}\right) $. Since $%
\gamma :=(f+\delta _{E})^{\ast }(x^{\ast })\leq \varphi _{\mathcal{H}%
}(x^{\ast })$ by (\ref{3.3k}), we have to prove that $\varphi _{\mathcal{H}%
}(x^{\ast })\leq \gamma $. This is obvious if $\gamma =+\infty $. Moreover,
since $E\cap \limfunc{dom}f\not=\emptyset $ we have $\gamma \not=-\infty $.
Suppose now that $\gamma \in \mathbb{R}$. Suppose now that $\gamma \in
\mathbb{R}$. We have $(x^{\ast },\gamma )\in \mathrm{epi}(f+\delta
_{E})^{\ast }$, and by Lemmas \ref{lem1} and \ref{lem51aa}$\left( \mathrm{ii}%
\right) $,
\begin{equation*}
(x^{\ast },\gamma )\in \left( \overline{\mathrm{co}}\mathcal{A}_{\mathcal{H}%
}\right) \cap \left( \{x^{\ast }\}\times \mathbb{R}\right) .
\end{equation*}%
Now, as $\left( \mathrm{ii}\right) $ holds, there exists a net $%
(r_{i})_{i\in I}$ such that
\begin{equation*}
\lim_{i\in I}r_{i}=\gamma ,\ \ (x^{\ast },r_{i})\in \mathcal{A}_{\mathcal{H}%
},\ \ \forall i\in I.
\end{equation*}%
Again, it follows from Lemma \ref{lem51aa}(i), that $\varphi _{\mathcal{H}%
}(x^{\ast })\leq r_{i}$ for all $i\in I$. Passing to the limit we get $%
\varphi _{\mathcal{H}}(x^{\ast })\leq \gamma $, as desired. \hfill $\square $

\begin{remark}
It is worth emphasizing that the implication $\left( \mathrm{i}\right)
\Rightarrow \left( \mathrm{ii}\right) $ of Theorem \ref{thm41a} holds for
arbitrary proper functions $f$ and$\ f_{t},$ $t\in T$.
\end{remark}

By Proposition \ref{pro1}, we have
\begin{equation*}
\mathcal{A}:=\bigcup\limits_{(\lambda _{t})\in \mathbb{R}_{+}^{(T)}}\mathrm{%
epi}\left( f+\sum\limits_{t\in T}\lambda _{t}f_{t}\right) ^{\ast }=\mathcal{A%
}_{\mathcal{F}(T)},
\end{equation*}%
and the corresponding $\mathcal{F}(T)$-dual problem of \textrm{(P)} reads:

\begin{equation*}
(\mathrm{D}\mathrm{)}\ \ \ \ \ \ \sup\limits_{ (\lambda_t) \in \mathbb{R}%
_{+}^{(T)}}\inf\limits_{x\in X}\left\{ f(x)+\sum\limits_{t\in T}\mu
_{t}f_{t}(x)\right\} . \qquad \qquad
\end{equation*}

\begin{corollary}
\label{corol61} Assume that $\{f;\ f_{t},\ t\in T\}\subset \Gamma (X)$ and $%
E\cap \limfunc{dom}f\not=\emptyset $. The following statements are
equivalent:\newline
$\left( \mathrm{i}\right) $ \ $\inf \mathrm{(P)}=\sup \mathrm{(D)}$. \newline
$\left( \mathrm{ii}\right) $ \ $\overline{\mathcal{A}}\cap \left(
\{0_{X^{\ast }}\}\times \mathbb{R}\right) =\overline{\mathcal{A}\cap \left(
\{0_{X^{\ast }}\}\times \mathbb{R}\right) }$.
\end{corollary}

\noindent \textbf{Proof } Applying Theorem \ref{thm41a} for $\mathcal{H=F}%
(T) $, one has $\mathcal{A}_{\mathcal{H}}=\mathcal{A},$ which is convex.
Consequently, $\overline{\mathrm{co}}\mathcal{A}_{\mathcal{H}}=\overline{%
\mathcal{A}}$ and, taking $x^{\ast }=0_{X^{\ast }}$, we are done. \hfill $%
\square $

We now come back to the general LIP problem in (\ref{LIP}), with $c^{\ast }$
instead of $x^{\ast },$

\begin{equation*}
(\mathrm{P})\ \ \ \ \ \ \inf \langle c^{\ast },x\rangle \ \ \text{ s.t.}\ \
\langle a_{t}^{\ast },x\rangle \leq b_{t},\ t\in T,
\end{equation*}%
and its $\mathcal{H}$-dual,%
\begin{equation*}
(\mathrm{D})\ \ \ \ \ \ \sup\limits_{\substack{ H\ \in \mathcal{H},\mu \in
\mathbb{R}_{+}^{H}\  \\ \sum\limits_{t\in H}\lambda _{t}a_{t}^{\ast
}=-c^{\ast }}}-\sum\limits_{t\in H}\lambda _{t}\,b_{t},
\end{equation*}

Recalling the sets $\mathcal{K}_{\mathcal{H}}$ and $\mathcal{A}_{\mathcal{H}%
} $ defined in (\ref{KH}) and (\ref{4.1bis}), respectively,
\begin{eqnarray*}
\mathcal{K}_{\mathcal{H}} &=&\tbigcup\limits_{H\in \mathcal{H}}\mathcal{K}%
_{H}=\tbigcup\limits_{H\in \mathcal{H}}\limfunc{cone}\left( \left\{
(a_{t}^{\ast },b_{t}),\ t\in H\right\} +\left\{ 0_{X^{\ast }}\right\} \times
\mathbb{R}_{+}\right) , \\
\mathcal{A}_{\mathcal{H}} &=&\left\{ \left( c^{\ast },0\right) \right\} +%
\mathcal{K}_{\mathcal{H}},
\end{eqnarray*}%
we now can state

\begin{corollary}
\label{corol42b} Assume that the LIP problem $\mathrm{(P)}$ is feasible and $%
\mathcal{H}$ is covering. Then the following statements are equivalent:%
\newline
$\mathrm{(i)}$ \ $\inf \mathrm{(P)}=\sup \mathrm{(D_{\mathcal{H}})}$.\newline
$\mathrm{(ii)}$ $\left( \overline{\mathrm{co}}\mathcal{K}_{\mathcal{H}%
}\right) \cap \left( \{-c^{\ast }\}\times \mathbb{R}\right) =\overline{%
\mathcal{K}_{\mathcal{H}}\cap \left( \{-c^{\ast }\}\times \mathbb{R}\right) }
$.
\end{corollary}

\noindent \textbf{Proof } Note that when $f\equiv 0$, $f_{t}=a_{t}^{\ast
}-b_{t}$ ($t\in T$), $x^{\ast }=-c^{\ast }$, we have
\begin{equation*}
-\inf (\mathrm{D})=\left( f+\delta _{E}\right) ^{\ast }(-c^{\ast }).
\end{equation*}%
The second part of Theorem \ref{thm41a} now concludes the proof. \hfill $%
\square \medskip $

In Example \ref{Exam4.2} with $x^{\ast }\in \mathbb{R}_{+}^{2}\diagdown
\left\{ 0_{2}\right\} ,$ by the characterization of $\mathcal{K}_{\mathcal{H}%
_{1}}$ in Example \ref{Exam4.3}, one has
\begin{equation*}
\begin{array}{ll}
\overline{\limfunc{co}}\mathcal{K}_{\mathcal{H}_{1}}\cap \left( \{-x^{\ast
}\}\times \mathbb{R}\right) & =\{-x^{\ast }\}\times \left[ -\frac{%
x_{1}^{\ast }x_{2}^{\ast }}{x_{1}^{\ast }+x_{2}^{\ast }},+\infty \right]
\medskip \\
& =\overline{\{-x^{\ast }\}\times \left] -\frac{x_{1}^{\ast }x_{2}^{\ast }}{%
x_{1}^{\ast }+x_{2}^{\ast }},+\infty \right[ }\medskip \\
& =\overline{\mathcal{K}_{\mathcal{H}_{1}}\cap \left( \{-x^{\ast }\}\times
\mathbb{R}\right) },%
\end{array}%
\end{equation*}%
so that condition $\mathrm{(ii)}$ in Corollary \ref{corol42b} holds. Thus, $%
\inf \mathrm{(P)}=\sup \mathrm{(D_{\mathcal{H}_{1}})}$. The above argument
is valid for Example \ref{Exam4.1} just replacing the interval $\left] -%
\frac{x_{1}^{\ast }x_{2}^{\ast }}{x_{1}^{\ast }+x_{2}^{\ast }},+\infty %
\right[ $ by its closure.

In the case where $\mathcal{H}$ is directed covering we have
\begin{equation*}
\mathcal{K}_{\mathcal{H}}=\mathcal{K}=\limfunc{cone}\left( \left\{
(a_{t}^{\ast },b_{t}),\ t\in T\right\} +\left\{ 0_{X^{\ast }}\right\} \times
\mathbb{R}_{+}\right) \text{.}
\end{equation*}%
The corresponding $\mathcal{F}(T)$-dual problems of (P) reads:

\begin{equation*}
(\mathrm{D})\ \ \ \ \ \ \sup\limits_{\substack{ (\lambda _{t})\ \in \
\mathbb{R}_{+}^{(T)}  \\ \sum\limits_{t\in T}\lambda _{t}a_{t}^{\ast
}=-c^{\ast }}}-\sum\limits_{t\in H}\lambda _{t}\,b_{t},
\end{equation*}%
which is the familiar Haar-dual problem of the LIP problem \textrm{(P)}.

\begin{corollary}
\label{corol43b} Assume that the LIP problem $(\mathrm{P})$ is feasible.
Then, the following statements are equivalent:\newline
$\mathrm{(i)}$ \ $\inf \mathrm{(P)}=\sup \mathrm{(D)}$.\newline
$\mathrm{(ii)}$ $\overline{\mathcal{K}}\cap \left( \{-c^{\ast }\}\times
\mathbb{R}\right) =\overline{\mathcal{K}\cap \left( \{-c^{\ast }\}\times
\mathbb{R}\right) }$.
\end{corollary}

\noindent \textbf{Proof } Apply Corollary \ref{corol42b} for $\mathcal{H}=%
\mathcal{F}(T)$. Then $\mathcal{K}_{\mathcal{H}}=\mathcal{K}$ is convex, $%
\overline{\limfunc{co}}\mathcal{K}=\overline{\mathcal{K}}$, and we are done.
\hfill $\square $

\begin{remark}
Corollary \ref{corol43b} has been quoted in \cite[Theorem 8.2]{GL98} for the
particular case when $(\mathrm{P})$ is an LSIP problem and $c^{\ast }$
belongs to the relative boundary of the $\limfunc{cone}\{a_{t}^{\ast },t\in
T\}$.
\end{remark}

\textbf{Funding}

This research was supported by the Vietnam National University HoChiMinh
city (VNU-HCM) under the grant number B2021-28-03, and by Ministerio de
Ciencia, Innovaci\'{o}n y Universidades (MCIU), Agencia Estatal de
Investigaci\'{o}n (AEI), and European Regional Development Fund (ERDF),
Project PGC2018-097960-B-C22.

\end{document}